\documentclass[12pt]{article} 
\usepackage{latexsym} 
\usepackage{epsfig} 
\newcommand{\sfrac}[2]{{\displaystyle{\frac{#1}{#2}}}}

\begin{document}

\title{Memoir on the Theory of the Articulated Octahedron }

\author{Author: Raoul Bricard(\cite{bio})}
\begin{small}
\date{Translator: Evangelos A. Coutsias, March 23, 2010 \footnote{
Translation from the French original of Raoul Bricard's 1897 memoir
on articulated octahedra \cite{Br}. 
(E. A. Coutsias, coutsias@unm.edu, Mathematics
Dept., University of New Mexico \copyright). Footnotes in the original
are listed in the bibliography. Notes by the translator are shown in 
parentheses. Section titles (in parentheses) were provided by the translator.
Translations of the replies to the problem of C. Stephanos by R. Bricard 
and C. Juel in the Interm\'ediaire
des Math\'ematiciens are also provided in the bibliography.
Only equations numbered are those numbered in the original. Changes in equation 
numbering are as follows:
 $(4,4^\prime)\rightarrow (4a,4b); (5,5^\prime) \rightarrow (5a,5b); 
(7^\prime,8^\prime,9^\prime) \rightarrow (12,13,14); 
(7^{\prime\prime},8^{\prime\prime},9^{\prime\prime}) \rightarrow (15,16,17); (12,13) \rightarrow (18,19)$.
}}

\end{small}

\maketitle

\vskip .10in
\section{(Introduction)}
Mr. C. Stephanos posed  the following question in the ``{\em Interm\'ediaire 
des Math\'ematiciens }'' \cite{IM1}:\\
``Do there exist polyhedra with invariant facets that are susceptible to an
infinite family of transformations that only alter solid angles and
dihedrals?''\\
I announced, in the same Journal \cite{IM2}, a special concave octahedron
possessing the required property. Cauchy, on the other hand, has proved
\cite{Cauchy} that there do not exist convex polyhedra that are deformable
under the prescribed conditions.\\
In this Memoir I propose to extend the above mentioned result, by resolving the
problem of Mr. Stephanos in general for octahedra of triangular facets.\\ 
Following Cauchy's theorem, all the octahedra which I shall establish as 
deformable will be of necessity concave by virtue of the fact that they possess 
reentrant dihedrals or, in fact, facets that intercross, in the manner of facets 
of polyhedra in higher dimensional spaces.
\section{(The tetrahedral equation and its decomposition)}
\begin{figure}
\centerline{\epsfig{figure=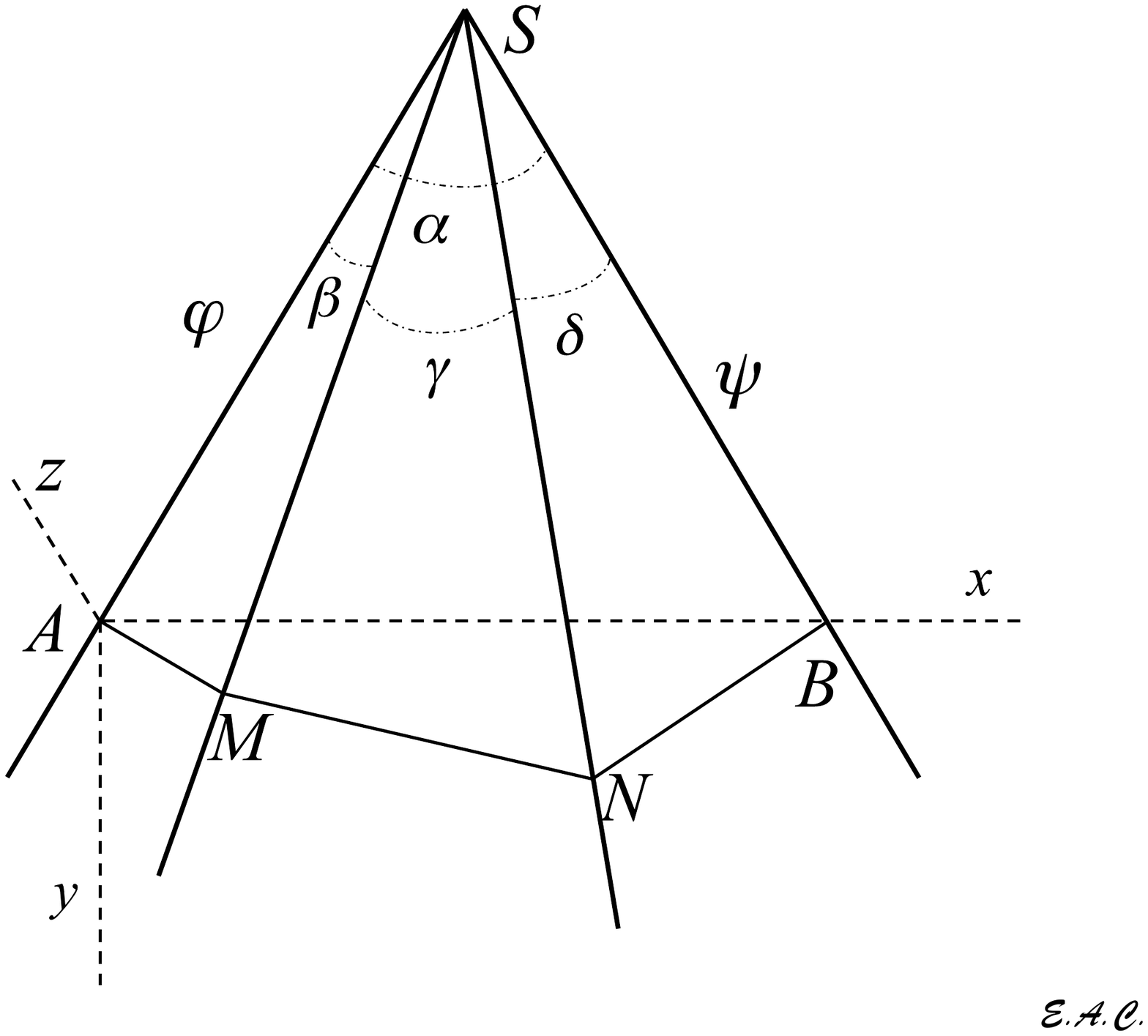,width=5.0in}}
\caption{}
\label{fig:tetrahedral}
\end{figure}

I shall begin by establishing certain properties pertaining to 
the deformation of
a tetrahedral angle, whose four faces remain invariant. This deformation is
analogous to that of a planar articulated quadrilateral. Given (fig.
\ref{fig:tetrahedral})
the tetrahedral $SABNM$, articulated along its four edges and having fixed-size
faces
\[ \angle ASB = \alpha\ ,\ \angle ASM = \beta\ ,\ \angle MSN = \gamma\ ,\ 
\angle NSB = \delta\ 
(0<\alpha\ ,\ \beta\ ,\ \gamma\ ,\ \delta<\pi) \ ,\]
let us find the relation between the dihedrals
$SA = \phi$ and $SB = \psi$.
One may suppose that the face $ASB$ maintains a fixed position. We associate
with it an orthogonal coordinate system of 3 axes defined as follows:
the origin is placed at a point $A$ of $SA$, so that $SA = 1$. The axis
$Ax$ ($Ay$) is parallel to the exterior (interior) bisectrix 
of the angle $\angle ASB$
and oriented so that the point $S$ has a negative ordinate and the point
$B$ a positive abscissa. The direction of positive $z$ may be taken arbitrarily.

With their axes thus chosen, the points $M$ and $N$, lying on the edges
$SM$, $SN$ and such that the angles $\angle SAM$, $\angle SBN$ are right 
angles, have as coordinates, respectively:
\[ {\bf M}
\left\{ \begin{array}{ccl}
x_1 & = & \tan{\beta}\cos{\frac{\alpha}{2}}\cos{\phi}\ , \\
y_1 & = & \tan{\beta}\sin{\frac{\alpha}{2}}\cos{\phi}\ , \\
z_1 & = & \tan{\beta}\sin{\phi}\ , \end{array}\right.
{\bf N}
\left\{ \begin{array}{ccl}
x_2 & = & 2\sin{\frac{\alpha}{2}} -
           \tan{\delta}\cos{\frac{\alpha}{2}}\cos{\psi}\ , \\
y_2 & = & \tan{\delta}\sin{\frac{\alpha}{2}}\cos{\psi}\ , \\
z_2 & = & \tan{\delta}\sin{\psi}\ . \end{array}\right.          \]
Equating the value of $\bar{MN}^2$ that results from these expressions
to that found by considering the triangle $\triangle SMN$, we get
\begin{eqnarray*}
\frac{1}{\cos^2{\beta}} + \frac{1}{\cos^2{\delta}} -
\frac{2\cos{\gamma}}{\cos{\beta}\cos{\delta}} & &\\
= \left(\tan{\beta}\cos{\frac{\alpha}{2}}\cos{\phi}
+       \tan{\delta}\cos{\frac{\alpha}{2}}\cos{\psi} - 2\sin{\frac{\alpha}{2}}
  \right)^2 & &\\
+\left(\tan{\beta}\sin{\frac{\alpha}{2}}\cos{\phi} 
-      \tan{\delta}\sin{\frac{\alpha}{2}}\cos{\psi}\right)^2 & &\\
+\left(\tan{\beta}\sin{\phi} - \tan{\delta}\sin{\psi}\right)^2& & \ ,
\end{eqnarray*}
and, after reductions,
\begin{eqnarray*}
  \sin{\beta}\sin{\delta}\cos{\alpha}\cos{\phi}\cos{\psi}
- \sin{\beta}\sin{\delta}\sin{\phi}\sin{\psi} & & \\
-\sin{\alpha}\sin{\beta}\cos{\delta}\cos{\phi} 
-\sin{\alpha}\sin{\delta}\cos{\beta}\cos{\psi} & & \\
+ \cos{\gamma} - \cos{\alpha}\cos{\beta}\cos{\delta} & = & 0 \ .
\end{eqnarray*}
Let us now introduce the variables 
\[ t = \tan{\frac{\phi}{2}} \ \ \mbox{and}\ \ u = \tan{\frac{\psi}{2}} \ .\]
We have
\begin{eqnarray*}
\cos{\phi} = \frac{1-t^2}{1+t^2}\ , & & \sin{\phi} = \frac{2t}{1+t^2}\ ,\\
\cos{\psi} = \frac{1-u^2}{1+u^2}\ , & & \sin{\psi} = \frac{2u}{1+u^2}\ .
\end{eqnarray*}
We arrive thus at the following relationship, which I shall term
``{\em equation of the tetrahedral angle}'':
\begin{equation}
\label{tetrahedral}
At^2u^2 + Bt^2 +2Ctu +Du^2 + E = 0 \ ,
\end{equation}
where
\begin{equation}
\left. \begin{array}{ccl}
A & = & \sin{\beta}\sin{\delta}\cos{\alpha} 
    +   \sin{\beta}\cos{\delta}\sin{\alpha} \\
  &   & +   \sin{\delta}\cos{\beta}\sin{\alpha} 
    -    \cos{\alpha}\cos{\beta}\cos{\delta} + \cos{\gamma} \\
  & = & \cos{\gamma} - \cos{\left(\alpha + \beta + \delta\right)}\ , \\
B & = & -\sin{\beta}\sin{\delta}\cos{\alpha}
    +   \sin{\beta}\cos{\delta}\sin{\alpha} \\
  &   & -   \sin{\delta}\cos{\beta}\sin{\alpha} 
    -    \cos{\alpha}\cos{\beta}\cos{\delta} + \cos{\gamma} \\
  & = & \cos{\gamma} - \cos{\left(\alpha + \beta - \delta\right)}\ , \\
C & = & -2\sin{\beta}\sin{\delta}\ , \\
D & = & -\sin{\beta}\sin{\delta}\cos{\alpha}
    -   \sin{\beta}\cos{\delta}\sin{\alpha} \\
  &   & +   \sin{\delta}\cos{\beta}\sin{\alpha} 
    -    \cos{\alpha}\cos{\beta}\cos{\delta} + \cos{\gamma} \\
  & = & \cos{\gamma} - \cos{\left(\alpha - \beta + \delta\right)}\ , \\
E & = & \sin{\beta}\sin{\delta}\cos{\alpha}
    -   \sin{\beta}\cos{\delta}\sin{\alpha} \\
  &   & -   \sin{\delta}\cos{\beta}\sin{\alpha} 
    -    \cos{\alpha}\cos{\beta}\cos{\delta} + \cos{\gamma} \\
  & = & \cos{\gamma} - \cos{\left(\alpha - \beta - \delta\right)}\ .
\end{array} \right\}
\end{equation}
It should be easy to anticipate the form of relation (\ref{tetrahedral}).

In fact, to a fixed value of $t = \tan{\left(\phi/2\right)}$ there
corresponds a unique position of the face $ASM$, the angle $\phi$ being
determined in this way up to an integral multiple of $2\pi$.
This face being thus fixed, the construction of the tetrahedral angle
can be achieved in two ways:
there exist, in fact, two positions of the line $SN$, symmetric with respect
to the plane $BSM$ and such that we have 
\[ \angle MSN = \gamma\ ,\ \angle BSN = \delta \ . \]
To every position of the face $BSN$ there corresponds a single value of the
variable $u$. Therefore the relationship that couples $t$ and $u$ must be of
the second degree with respect to $u$. 

For the same reason, that relationship should be of the second degree with
respect to $t$. At last it can be seen that if this relation is satisfied 
by the values $t$, $u$, it is equally true for $-t$ and $-u$. It is therefore
necessarily of the form (\ref{tetrahedral}).

{\em Cases of decomposition of equation (\ref{tetrahedral})}.-- We just saw
that to a value of $t$ there correspond two values of $u$. It is useful to 
investigate under what conditions relation (\ref{tetrahedral}) decomposes,
in a way that these values of $u$ are rational functions of $t$.

For that to be so, it is necessary and sufficient that the polynomial
\[ C^2t^2-\left(At^2+D\right)\left(Bt^2+E\right) =
-ABt^4 + \left(C^2-AE-BD\right)t^2 -DE \ , \]
that appears under the radical entering in the expression for $u$ as a
function of $t$, is a perfect square. This can be achieved in two ways:
\begin{enumerate}
\item There holds
\[ \left(C^2-AE-BD\right)^2-4ABDE = 0 \ . \]
It is found, by an easy calculation, that the left hand side of this equality
reduces to the expression
\[ 16 \sin^2{\alpha}\sin^2{\beta}\sin^2{\gamma}\sin^2{\delta} \ .\]
This condition is not satisfied unless one of the angles $\alpha$, $\beta$,
$\gamma$ or $\delta$ is reduced to $0$ or $\pi$, which is impossible.
It would seem that this case presents itself when the vertex $S$ is moved to
infinity and the tetrahedral degenerates to a prismatic solid. But then one
can see that the values of the coefficients $A$, $B$, $C$, $D$, $E$, are all
reduced to $0$, and the preceding reasoning no longer applies. Any right
section of the prismatic solid is an articulated quadrilateral, whose 
deformation is governed by an equation of the same form as equation
(\ref{tetrahedral}):
\[
A't^2u^2 + B't^2 +2C'tu +D'u^2 + E' = 0 \ .
\]
The coefficients $A'$, $B'$, $C'$, $D'$, $E'$ have, in terms of the edges
$a$, $b$, $c$, $d$ of the quadrilateral, expressions which can be obtained
from the forms of $A$, $B$, $C$, $D$, $E$ by allowing in the latter the
angles $\alpha$, $\beta$, $\gamma$, $\delta$ to approach zero in such a way
that 
\[ \frac{\alpha}{a} = \frac{\beta}{b} = \frac{\gamma}{c} = \frac{\delta}{d}\ .
\]
One can then see that the condition
\[ \left(C'^2-A'E'-B'D'\right)^2-4A'B'D'E' = 0 \ , \]
implies the impossible relationship
\[ abcd = 0 \ .\]
\item There holds
\[ AB = 0\ \  \mbox{    with    }\ \  DE = 0 \ . \]
These relationships imply one of the following cases:
\begin{equation}
\label{caseII}
\left. \begin{array}{ll}
B = 0\ , & E = 0\ ,\\
A = 0\ , & D = 0\ ,\\
B = 0\ , & D = 0\ ,\\
A = 0\ , & E = 0\ . \end{array} \right\}
\end{equation}
Let us consider for example the equalities
\[ A = 0\ , \ \ \ \ \ \ \ \  D = 0 \ . \]
There results:
\begin{eqnarray*}
\mbox{let  }\ 
\gamma = \alpha+\beta+\delta+2k\pi & \mbox{with} &
\gamma = \alpha-\beta+\delta+2k'\pi \\
\mbox{let  }\ 
\gamma = \alpha+\beta+\delta+2k\pi & \mbox{with} &
\gamma = -\alpha+\beta-\delta+2k'\pi \\
\mbox{let  }\ 
\gamma = -\alpha-\beta-\delta+2k\pi & \mbox{with} &
\gamma = \alpha-\beta+\delta+2k'\pi \\
\mbox{let  }\ 
\gamma =-\alpha-\beta-\delta+2k\pi & \mbox{with} &
\gamma =-\alpha+\beta-\delta+2k'\pi \ .
\end{eqnarray*}
The first pair of relationships is incompatible with the hypotheses on
the magnitudes of the angles $\alpha$, $\beta$, $\gamma$, $\delta$.
It can be seen, in fact, that
\[ \beta +\left(k-k'\right)\pi = 0 \ ,\]
which is impossible.

In examining the other pairs, it is seen that only the third is
admissible and that it has as necessary consequences
\[ \alpha + \delta = \pi\ ,\ \ \ \ \beta+\gamma = \pi \ .\]
It will be argued similarly on the remaining equalities
(\ref{caseII}). I will rewrite them anew, listing next to each the relationship
implied among the faces of the tetrahedral:
\[
\begin{array}{llll}
B = 0\ ,& E = 0\ ,& \delta = \alpha\ , & \gamma = \beta \\
A = 0\ ,& D = 0\ ,& \delta = \pi - \alpha\ , & \gamma = \pi - \beta \\
B = 0\ ,& D = 0\ ,& \delta = \beta\ , & \gamma = \alpha \\
A = 0\ ,& E = 0\ ,& \delta = \pi - \beta\ , & \gamma = \pi - \alpha \ .
\end{array}
\]
\end{enumerate}
We are now able to recognize two cases of decomposing equation
(\ref{tetrahedral}):
\begin{enumerate}
\item The tetrahedral has its adjacent faces equal or supplementary in
pairs. Its equation reduces to
\[ At^2u + 2Ct + Du = 0\ ,\]
or
\[ Bt^2 + 2Ctu + E = 0\ .\]
(by omitting in the first the factor $u=0$ which corresponds to the
uninteresting case where the adjacent faces of the tetrahedral remain
coincident during the deformation). I will say that a tetrahedral of this
nature is {\em rhomboidal}.
\item The tetrahedral has opposing angles equal or supplementary in pairs.
Its equation is then
\[ At^2u^2 + 2Ctu + E = 0\ ,\]
or
\[ Bt^2 + 2Ctu + Du^2 = 0\ .\]
\end{enumerate}
These equations are written, respectively, by introducing the explicit 
forms of their coefficients,
\[
\left[\cos{\alpha} -
      \cos{\left(\alpha+2\beta\right)}\right]t^2u^2
-4\sin^2{\beta}tu+\cos{\alpha}-\cos{\left(\alpha-2\beta\right)} =0 \ ,\]
or
\[ \sin{\left(\alpha+\beta\right)}t^2u^2-2\sin{\beta}tu-\sin{\left(\alpha-
\beta\right)} = 0 \]
and
\[
\left[\cos{\left(\alpha+2\beta\right)} - \cos{\alpha}\right]t^2
-2\sin^2{\beta}tu+ \left[\cos{\left(\alpha-2\beta\right)} - \cos{\alpha}\right]
u^2 = 0 \ ,\]
or
\[
\sin{\left(\alpha+\beta\right)}t^2 + 2\sin{\beta}tu -\sin{\left(\alpha-
\beta\right)}u^2 = 0 \ .\]
They decompose, the first into
\begin{eqnarray}
& & tu =  \frac{\sin{\beta}+\sin{\alpha}}{\sin{\left(\alpha+\beta\right)}}
   = \frac{\cos{\frac{\alpha-\beta}{2}}}{\cos{\frac{\alpha+\beta}{2}}} 
\ \ \ \ \ \ \left(a\right) \nonumber \\
 \mbox{and} & & \\
& & tu = \frac{\sin{\beta}-\sin{\alpha}}{\sin{\left(\alpha+\beta\right)}}
   = \frac{\sin{\frac{\beta-\alpha}{2}}}{\sin{\frac{\alpha+\beta}{2}}}
\ \ \ \ \ \ \left(b\right) \nonumber 
\end{eqnarray}
the second into
\begin{eqnarray}
& & 
\frac{t}{u} = \frac{-\sin{\beta}+\sin{\alpha}}{\sin{\left(\alpha+\beta\right)}}
   = \frac{\sin{\frac{\alpha-\beta}{2}}}{\sin{\frac{\alpha+\beta}{2}}}
\ \ \ \ \ \ \left(a\right) \nonumber \\
\mbox{and} & & \\
& & 
\frac{t}{u} = \frac{-\sin{\beta}-\sin{\alpha}}{\sin{\left(\alpha+\beta\right)}}
   =-\frac{\cos{\frac{\alpha-\beta}{2}}}{\cos{\frac{\alpha+\beta}{2}}}
\ \ \ \ \ \ \left(b\right) \nonumber 
\end{eqnarray}
It is not unhelpful to summarize the previous discussion: we can distinguish
three types of articulated tetrahedral angles:
\begin{enumerate}
\item The {\em general} tetrahedral angles whose faces have no special
relation. Its equation is irreducible, of the form that to each value
of one of the variables $t$, $u$ there correspond two values of the other
variable, that are not rational functions of the first;
\item The {\em rhomboidal} tetrahedral angles. To one value of $t$ there
corresponds a single value of $u$, which is a rational function of $t$,
but the converse is not true;
\item The tetrahedral angles whose opposite faces are {\em equal or
supplementary in pairs}. To one value of $t$ there corresponds a unique
value of $u$ and conversely (\cite{uni}).
\end{enumerate}

\section{(Reconstruction and equivalence)}

Since equation (\ref{tetrahedral}) involves four arbitrary parameters,
any equation of similar form
\[
At^2u^2 + Bt^2 +2Ctu +Du^2 + E = 0 \ ,
\]
may be considered to define the deformation of an articulated tetrahedral
angle.

The elements of this tetrahedral angle are given by the relationships
\begin{eqnarray*}
\frac{\cos{\gamma} - \cos{\left(\alpha + \beta + \delta\right)}}{A} & = &
\frac{\cos{\gamma} - \cos{\left(\alpha + \beta - \delta\right)}}{B} =
\frac{-2\sin{\beta}\sin{\delta}}{C} \\
& = &
\frac{\cos{\gamma} - \cos{\left(\alpha - \beta + \delta\right)}}{D} =
\frac{\cos{\gamma} - \cos{\left(\alpha - \beta - \delta\right)}}{E}
\end{eqnarray*}
from which we can show
\begin{eqnarray*}
\frac{-2\sin{\beta}\sin{\delta}}{C} & = & 
\frac{4\sin{\beta}\sin{\delta}\cos{\alpha}}{A-B-D+E} =
\frac{4\sin{\delta}\cos{\beta}\sin{\alpha}}{A-B+D-E} \\
& = &
\frac{4\sin{\beta}\cos{\delta}\sin{\alpha}}{A+B-D-E} =
\frac{4\left(\cos{\gamma} - \cos{\alpha}\cos{\beta}\cos{\delta}\right)}
     {A+B+D+E} \ .
\end{eqnarray*}
We have, as a result:
\begin{equation}
\label{angles}
\left.
\begin{array}{c}
\cos{\alpha}=-\sfrac{A-B-D+E}{2C}\ ,\ \ 
\tan{\beta} = \sfrac{-2C\sin{\alpha}}{A-B+D-E}\ ,\\
\tan{\delta} = \sfrac{-2C\sin{\alpha}}{A+B-D-E}\ ,\\
\cos{\gamma} = \cos{\alpha}\cos{\beta}\cos{\delta} -
\sfrac{A+B+D+E}{2C} \sin{\beta}\sin{\delta}\ ,
\end{array}
\right\}
\end{equation}
formulas which permit the successive calculation of the angles
$\alpha$, $\beta$, $\delta$, $\gamma$.
Clearly, certain reality conditions need to be satisfied; however they
are quite complicated and it will not be of interest to state them here.

Since it can be assumed that
\[ \alpha\ ,\ \beta\ ,\ \gamma\ ,\ \delta\ <\ \pi \ ,
\]
the previous formulas define only {\em two} systems of values for the angles
(neglecting the second case of decomposition to which I shall return shortly).
It is seen immediately that if one of the systems
is formed from the values
\[ \alpha\ ,\ \beta\ ,\ \gamma\ ,\ \delta\ \ ,
\]
those of the other system are
\[ \alpha\ ,\ \pi - \beta\ ,\ \gamma\ ,\ \pi - \delta\ \ .
\]
There results thus a theorem which we will find extremely useful in the
sequel:\\
{\em If two articulated tetrahedral angles ${\bf T}$ and ${\bf T}_1$ can be 
subjected to a continuous
deformation in such a manner that two adjacent dihedrals of ${\bf T}$
remain equal or supplementary 
to two adjacent dihedrals of ${\bf T}_1$, these two tetrahedral
angles have all their faces pairwise equal or supplementary.}

Let us reserve the previous symbols for the elements of the tetrahedral
angle ${\bf T}$ and let us designate the corresponding elements of ${\bf T}_1$
by the same subscripted letters. The theorem posits three cases, 
all of which are established in similar fashion.

Let us suppose, for example, that we have invariably
\[ \phi = \phi_1\ ,\ \psi = \pi - \psi_1 \ ,\]
from which
\[ t_1 = t \ ,\  u_1 = \frac{1}{u} \ . \]
We have, during the deformation of ${\bf T}_1$, the relation
\[
A_1t_1^2u_1^2 + B_1t_1^2 +2C_1t_1u_1 +D_1u_1^2 + E_1 = 0 \ .
\]
Upon substituting $t$ and $\sfrac{1}{u}$ for $t_1$ and $u_1$, respectively,
this becomes  
\[
A_1\frac{t^2}{u^2} + B_1t^2 +2C_1\frac{t}{u} +D_1\frac{1}{u^2} + E_1 = 0 \ ,
\]
or
\[
B_1t^2u^2 + A_1t^2 +2C_1tu +E_1u^2 + D_1 = 0 \ .
\]
This relationship must be identical to (\ref{tetrahedral}). We have therefore
\[ \frac{B_1}{A} = \frac{A_1}{B} = \frac{C_1}{C} = \frac{E_1}{D} = 
\frac{D_1}{E} \ .
\]
Let us now apply the formulas in (\ref{angles}) to the tetrahedral angle
${\bf T}_1$. We find
\[
\alpha_1 = \pi - \alpha\ ,\ \beta_1 = \beta\ \mbox{or}\ \pi - \beta\ ,\ 
\delta_1 = \pi - \delta\ \mbox{or}\ \delta\ ,\ \gamma_1 = \gamma\ ,
\]
in agreement with the theorem stated above.

This proposition is still true when the tetrahedral angles
${\bf T}$ and ${\bf T}_1$ are rhomboidal, but it ceases to be if they
have each of their opposite faces equal or supplementary in pairs.
There exists in fact an infinity of such tetrahedral angles whose 
deformation is governed by the same equation
\[ tu = k \ \mbox{or}\ \frac{t}{u} = k^\prime\ . \]
If $\alpha$ and $\beta$ denote two adjacent faces of a tetrahedral angle
which satisfy the first relationship, e.g., we must have
\[
\sfrac{\cos{\frac{\beta-\alpha}{2}}}{\cos{\frac{\beta+\alpha}{2}}} = k\ 
\mbox{or, alternatively}\ \ 
\sfrac{\sin{\frac{\beta-\alpha}{2}}}{\sin{\frac{\beta+\alpha}{2}}} = k
\]
from which
\[
\tan{\frac{\alpha}{2}}\tan{\frac{\beta}{2}} = \frac{k-1}{k+1}\ 
\mbox{or, alternatively}\ \ 
\sfrac{\tan{\frac{\alpha}{2}}}{\tan{\frac{\beta}{2}}}
= \sfrac{1-k}{k+1}\ , \]
equalities which are satisfied for an infinity of values of $\alpha$ and 
$\beta$.
\section{(Covariance of opposite dihedrals in a tetrahedral)}
In order to complete these generalities, I shall establish now the following
property of the articulated tetrahedral angle, analogous to a well known 
property of the articulated quadrilateral.

{\em When an articulated tetrahedral angle is deformed, there exists a linear
relationship between the cosines of two opposite dihedrals}.

In effect, maintaining the previous notation and designating additionally
by $\theta$ the dihedral $ON$, we have, by the fundamental formulas of
spherical Trigonometry,
\begin{eqnarray*}
 \cos{BOM} & = & \cos{\alpha}\cos{\beta}
+\sin{\alpha}\sin{\beta}\cos{\phi} \\
 & = & \cos{\gamma}\cos{\delta}+\sin{\gamma}\sin{\delta}\cos{\theta}\ ,
\end{eqnarray*}
that is, a relation of the form
\[ A\cos{\phi} + B\cos{\theta} + C = 0 \ .\]
When the tetrahedral angle fits any of the cases of decomposition
listed above, this relationship reduces to
\[ \cos{\phi} = \cos{\theta} \, \]
from which
\[ \phi = \theta\ ,\]
if we do not allow for the dihedrals of the tetrahedral angle other than
values between $0$ and $\pi$. Conversely, if an articulated tetrahedral angle
is deformed in such a fashion that two opposite dihedrals remain always
equal, this tetrahedral angle is rhomboidal or it has its opposite faces 
equal or supplementary pairwise. In fact, we have now the relations
\[ \cos{\alpha}\cos{\beta} = \cos{\gamma}\cos{\delta}\ , \] 
\[ \sin{\alpha}\sin{\beta} = \sin{\gamma}\sin{\delta}\ , \] 
from which
\[ \cos{\left(\alpha\pm\beta\right)} = \cos{\left(\gamma\pm\delta\right)}\ ,\]
which admit the four sets of solutions
\begin{eqnarray*}
\gamma=\alpha\ ,\ \ \delta=\beta\ ;& & \gamma=\beta\ ,\ \ \delta=\alpha\ ;\\
\gamma=\pi-\alpha\ ,\ \ \delta=\pi-\beta\ ;& & \gamma=\pi-\beta\ ,
\ \ \delta=\pi-\alpha\ .
\end{eqnarray*}

\section{(Deformability in general)}
We shall now undertake the study of the deformability conditions for an
octahedron with triangular facets, with edges of fixed lengths.

It must first be noted that such an octahedron, although concave, is
generally rigid. This follows from Legendre's theorem by which the number
of conditions necessary for the determination of a polyhedron are exactly
equal to the number of its edges. In effect, the proof of that theorem relies
entirely on the fact that the polyhedron is determined by the relation
of Euler (or of Descartes) and does not depend on its convexity or concavity.
That is, an octahedron with triangular facets is well defined by that relation,
whatever the disposition of its facets may be.

An octahedron whose edges are given is therefore in general completely
determined, and therefore undeformable. Our approach is to examine if in certain
cases, by means of particular relations existing between its edges, that
determinacy ceases to hold. Then the octahedron will be deformable.

\begin{figure}
\centerline{\epsfig{figure=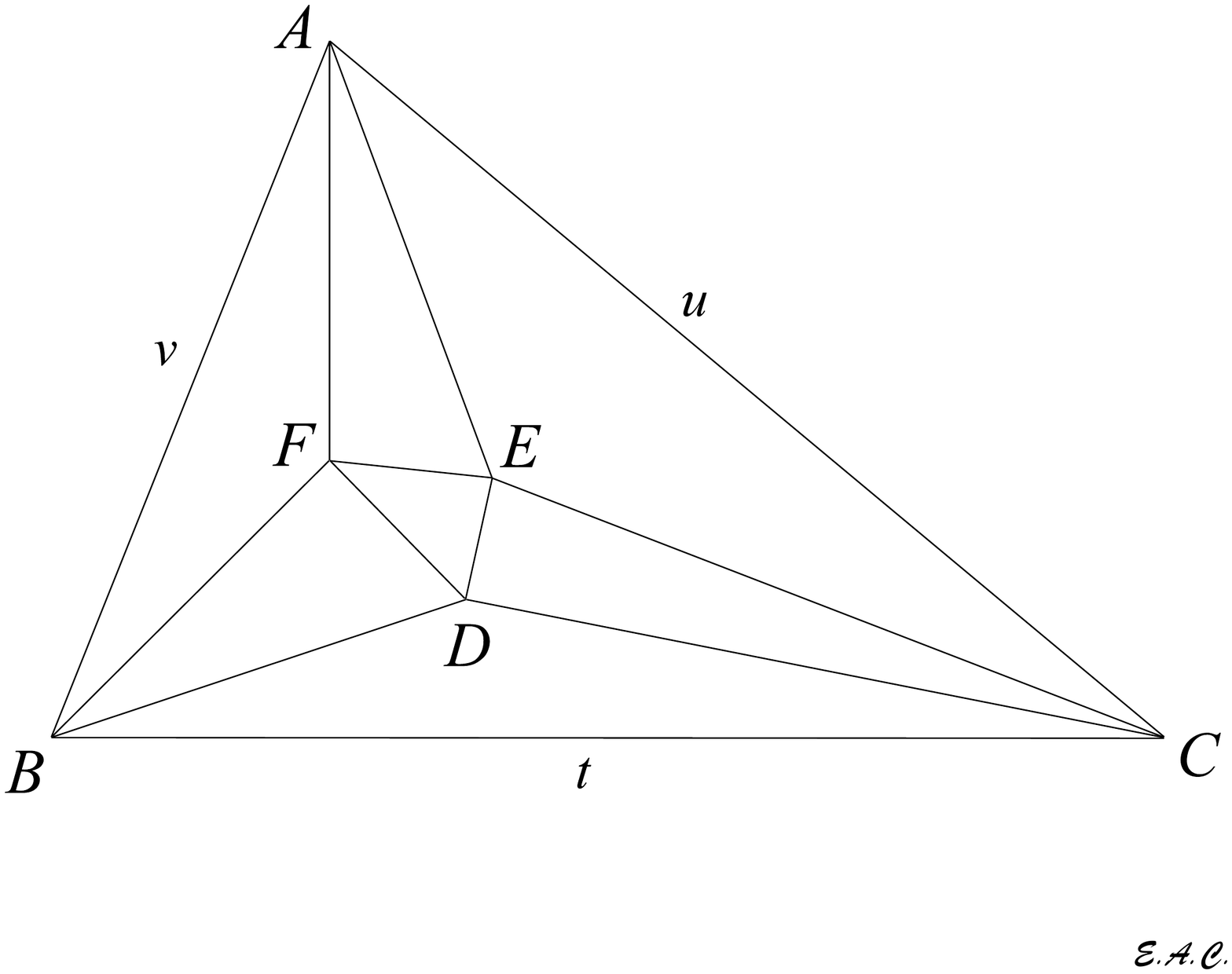,width=4.4in}}
\caption{}
\label{fig:triangle1}
\end{figure}

Let us suppose that this happens for the octahedron $ABCDEF$ 
(fig. \ref{fig:triangle1}). We can see then that the twelve dihedrals of this
octahedron are of necessity all variable, when it is subjected to the 
deformation to which it is amenable.

Let us assume in fact that, during the deformation of the octahedron
one of its dihedrals, for example $AB$, remains constant in magnitude.
The tetrahedral angle formed by the four faces having the point $A$
as common vertex will be completely rigid, since one of its dihedrals
is invariant. Consideration of the tetrahedral angles having their vertices
at $F$, $E$, $D$ shows then that all the dihedrals of the octahedron have
constant magnitude, which is contrary to the hypothesis.

The octahedron therefore is comprised of six tetrahedral angles which all
deform while keeping their faces constant. Three cases must be distinguished,
according to whether these tetrahedral angles are {\em general} (in the sense
given in \S II of this paper), {\em rhomboidal} or {\em having  opposite faces 
pairwise equal or supplementary}.

\section{(The case of a general tetrahedral)}
We examine now the first case. Among the six tetrahedral angles, let us 
consider those whose vertices are at $A$, $B$, $C$. Their deformations
are governed by three equations that are similar to (\ref{tetrahedral}) and
{\it indecomposable}:

\begin{equation}
\label{tetra}
At^2u^2 + Bt^2 +2Ctu +Du^2 + E = 0 \ ,
\end{equation}
\begin{equation}
\label{tetrb}
A^\prime t^2v^2 + B^\prime t^2 +2C^\prime tv +D^\prime v^2 + E^\prime  = 0 \ ,
\end{equation}
\begin{equation}
\label{tetrc}
A^{\prime\prime} u^2v^2 + B^{\prime\prime} u^2 +2C^{\prime\prime} uv 
+D^{\prime\prime} v^2 + E^{\prime\prime}  = 0 \ ,
\end{equation}

\noindent by designating as $t$, $u$, $v$ the tangents of the semi-dihedrals
$BC$, $CA$, $AB$ and as $A$, $B$, $\ldots$, $E''$ constants that depend on the 
faces of the three tetrahedral angles and, consequently, on the edges of the
octahedron.

The preceding equations ought to be satisfied by an infinity of sets of values
of $t$, $u$ and $v$. Therefore, equations (\ref{tetrb}) and (\ref{tetrc}),
in terms of $v$, ought to have {\it one} or {\it two} common roots for an
infinity of values of $t$ and of $u$ that satisfy equation (\ref{tetra}).

However it is impossible that equations (\ref{tetrb}) and (\ref{tetrc})
have their two roots always equal: indeed if it was so, we would have
\[ \frac{A't^2+D'}{A''u^2+D''} = \frac{C't}{C''u} =
\frac{B't^2+E'}{B''u^2+E''} \ , \]
from which one obtains two equations of third degree in $t$ and $u$, which 
ought to have, with equation (\ref{tetra}), an infinity of common solutions;
however this is impossible, since the latter was supposed indecomposable.

Thus equations (\ref{tetrb}) and (\ref{tetrc}) have, in general, a unique
common root in $v$. This root is therefore expressible as a rational function
of the coefficients of these equations and, consequently, of $t$ and $u$.
Now, one can derive from equations (\ref{tetra}) and (\ref{tetrb})
\begin{eqnarray*}
u & = & \frac{-Ct\pm \sqrt{F(t)}}{At^2+D}\ , \\
v & = & \frac{-C't\pm \sqrt{F_1(t)}}{A't^2+D'}\ ,
\end{eqnarray*}
defining
\begin{eqnarray*} 
F(t)   & = & C^2t^2-\left(At^2+D\right)\left(Bt^2+E\right)\ ,\\
F_1(t) & = & C'^2t^2-\left(A't^2+D'\right)\left(B't^2+E'\right)\ ,
\end{eqnarray*}
and conveniently choosing the signs placed in front of the radicals.
One therefore has
\[ \frac{-C't\pm \sqrt{F_1(t)}}{A't^2+D'}
= \phi\left[t,\frac{-Ct\pm \sqrt{F(t)}}{At^2+D}\right] \ ,
\]
where $\phi(x,y)$ denotes a rational function (\cite{typo}).

The second member of this relationship may be put in the form
\[ \frac{L+M\sqrt{F(t)}}{N}\ ,\]
$L$, $M$, $N$ being polynomials in $t$. One arrives finally at the 
identity
\[ P\sqrt{F(t)} + Q \sqrt{F_1(t)} + R = 0\ ,\]
where $P$, $Q$, $R$ are also polynomials in $t$. One extracts from that
\[ F(t)F_1(t) = \left[\frac{R^2-P^2F(t)-Q^2F_1(t)}{2PQ}\right]^2\ .
\]
The product of the polynomials $F(t)$ and $F_1(t)$ must therefore be the
square of a rational function and, consequently, of a polynomial in $t$.
It there follows that $F(t)$ and $F_1(t)$ are identical up to a constant
factor.

Indeed $F(t)$ and $F_1(t)$ are two biquadratic polynomials that are not
perfect squares, else equations (\ref{tetra}) and (\ref{tetrb}) would be
reducible, contrary to the stated hypothesis. One may therefore set
\begin{eqnarray*}
F(t) & = & -A\ B\ (t-\lambda\,)(t+\lambda\,)(t-\mu\,)(t+\mu\,)\ ,\qquad
\lambda\, \ne \mu\, ,\\
F_1(t) & = & -A'B'(t-\lambda')(t+\lambda')(t-\mu')(t+\mu')\ ,\qquad
\lambda' \ne \mu' ,
\end{eqnarray*}
and their product cannot be a perfect square unless one has
\[ \lambda=\pm\lambda\,',\qquad \mu = \pm\mu\,',\]
or 
\[ \lambda=\pm\mu\,',\qquad \mu = \pm\lambda\,',\]
which establish the above assertion.
We deduce from this the following important consequence:\\
\indent The equations (\ref{tetra}) and (\ref{tetrb}), respectively in terms
of $u$ and $v$, have their roots equal for the same values of $t$.

I could pursue the algebraic study of the system (\ref{tetra}), (\ref{tetrb}),
(\ref{tetrc}), whose consideration ought to be by itself sufficient, as one
can see easily, to give the conditions of deformability of the octahedron.
But this would lead to rather involved calculations, by reason of the
complicated dependence of the coefficients $A$, $B$, $C,\ldots$ on the 
elements of the octahedron. Therefore I shall take a different route.
I shall state however the following theorem, since it may find application
in other investigations:

{\em For equations (\ref{tetra}), (\ref{tetrb}),
(\ref{tetrc}), to have an infinity of common solutions, it is necessary and
sufficient that they result from the successive elimination of 
$t$, $u$, $v$ between the $2$ equations
\begin{eqnarray*}
luv+mvt+ntu+p & = & 0,\\
l't+m'u+n'v+p'tuv & = & 0,
\end{eqnarray*}
where $l,m,n,p,l',m',n',p'$ are arbitrary coefficients}.

Let us return then to the consideration of the octahedron $ABCDEF$ and let us
interpret geometrically the last result that was obtained.\\
\indent When $t$ assumes such a value that equation (\ref{tetra}) in $u$ 
has its roots equal, the dihedral $CE$ will become evidently equal to $0$
or $\pi$. Similarly, when equation (\ref{tetrb}) in $v$ has its roots equal,
the dihedral $BF$ will become equal to $0$ or $\pi$. Therefore the 
dihedrals $CE$ and $BF$ are such that if one of them becomes $0$ or $\pi$,
the other assumes one of these values simultaneously.\\
\indent Now, during the deformation of the octahedron, there exists a linear
relationship between the cosines of these two dihedrals. On has in effect
(IV) a linear relationship between the cosine of each of these
dihedrals and that of the dihedral $BC$, which opposes them in the articulated
tetrahedral angles with vertices respectively at $C$ and $B$ 
\begin{eqnarray*}
l\,\cos{CE} + m\,\cos{BC}+n\, & = & 0\,,\\
l'\cos{BF} + m'\cos{BC}+n' & = & 0\,
\end{eqnarray*}
from which there arises one more relation
\begin{equation}
\label{eq10}
l^{\prime\prime}\cos{CE} + m^{\prime\prime}\cos{BF}+n^{\prime\prime}  =  0\,
\end{equation}
with constant coefficients $l^{\prime\prime}\,,m^{\prime\prime}\,,
n^{\prime\prime}$.\\
\noindent Let us successively set in this relation the dihedral $CE$ equal
to $0$ or $\pi$. The dihedral $BF$, as we have stated, will assume each time
one of these values; it cannot assume the value $0$ or $\pi$ twice, because 
one must also have
\begin{eqnarray*}
 l^{\prime\prime}\pm m^{\prime\prime}+n^{\prime\prime} & = & 0\,,\\
-l^{\prime\prime}\pm m^{\prime\prime}+n^{\prime\prime} & = & 0\,,
\end{eqnarray*}
$m^{\prime\prime}$ having the same sign in each right side and, as a result,
\[ l^{\prime\prime} = 0\,,\qquad m^{\prime\prime} = \pm n^{\prime\prime}\,.\]
Relationship (\ref{eq10}) will reduce then to
\[ \cos{BF} = \pm 1\,, \]
which is impossible, since all the dihedrals of the octahedron are
variable. The dihedral $BF$ must therefore assume once the value $0$
and once the value $\pi$, in the same order as the dihedral $CE$ or in the 
reverse order. One has then
\begin{eqnarray*}
 l^{\prime\prime}\pm m^{\prime\prime}+n^{\prime\prime} & = & 0\,,\\
-l^{\prime\prime}\mp m^{\prime\prime}+n^{\prime\prime} & = & 0\,,
\end{eqnarray*}
with correspondence of the signs in the two right hand sides. One obtains from
this
\[ n^{\prime\prime} = 0\,,\qquad l^{\prime\prime} = \pm m^{\prime\prime}\,,\]
and relation (\ref{eq10}) becomes
\[ \cos{CE} = \pm \cos{BF}\,. \]
\indent Thus, during the deformation of the octahedron, {\em the dihedrals
$CE$ and $BF$ are invariably equal or supplementary}.\\
\indent We may employ the same reasoning in considering three tetrahedral
angles having as vertices those of a facet other than $ABC$. Since all these
tetrahedrals were assumed general, the conclusion will remain the same, and
we may state:

{\em During the deformation of the octahedron, its opposing dihedrals
will remain equal or supplementary in pairs}.

Let us then envision two opposing tetrahedral angles, for example those having
their vertices respectively at $A$ and $D$. These can be deformed in such a 
manner that their dihedrals remain invariably equal or supplementary in
pairs. Following the theorem of Sec. III, this implies that their faces must
be pairwise equal or supplementary. It is likewise for the faces of the
other three couples of opposed tetrahedral angles that constitute the
octahedron.\\
\indent It is easy to see that, solely, the equality of corresponding
faces is admissible. Let us consider, in effect, two opposing facets of the
octahedron, $\triangle ABC$ and $\triangle DEF$ for example. 
One has, as a result of the preceding discussion,
\begin{eqnarray*} 
\angle A = \angle D &\mbox{or} & \angle A+\angle D = \pi\,,\\
\angle B = \angle E &\mbox{or} & \angle B+\angle E = \pi\,,\\
\angle C = \angle F &\mbox{or} & \angle C+\angle F = \pi\,,
\end{eqnarray*} 
and, as it is shown in elementary geometry in order to establish the
similarity of two triangles having their sides pairwise parallel or 
perpendicular, the equalities written at the begining of each line 
are the only ones that may hold true.\\
\indent The octahedron is thus such that its opposing facets are
similar triangles, with homologous sides having always opposing edges.
One has, as a result, the sequence of equalities
\begin{eqnarray*} 
\frac{AB}{DE}=\frac{BC}{EF}=\frac{CA}{FD}\,,&\qquad &
\frac{CA}{FD}=\frac{AE}{DB}=\frac{EC}{BF}\,,\\
\frac{EC}{BF}=\frac{CD}{FA}=\frac{DE}{AB}\,,&\qquad &
\frac{CD}{FA}=\frac{DB}{AE}=\frac{BC}{EF}\,,
\end{eqnarray*} 
from which there follows
\begin{equation}
\label{typeI}
\left\{
\begin{array}{lll}
AB=DE\,, & BC = EF\,, & CA=FD\,, \\
AE=BD\,, & BF = CE\,, & CD=AF\,.
\end{array}\right.
\end{equation}

{\em Therefore the octahedron has its opposing edges pairwise equal}.

 I shall now show that these relationships suffice to ensure the deformability
of the octahedron (which does not follow from the preceding analysis),
{\em if one assumes additionally certain conditions related to the relative
placement of the facets}. There exist in fact, convex octahedra for which
the opposing edges are pairwise equal, and which, by Cauchy's theorem, cannot be
deformable.
\section{(Flexibility of Type I: axis of rotational symmetry)}
Let us consider to this effect a system of four invariant triangles
$\triangle AFB$, $\triangle DFB$, $\triangle ACE$, $\triangle DCE$ 
({\em fig. $3$}), joined at the points $A$ and $D$
and along the lines $BF$ and $CE$. One may suppose
\[ AF=DC\,, AB=DE\,, DB=AE\,, DF=AC\,, BF=CE\ . \]
\begin{figure}
\centerline{\epsfig{figure=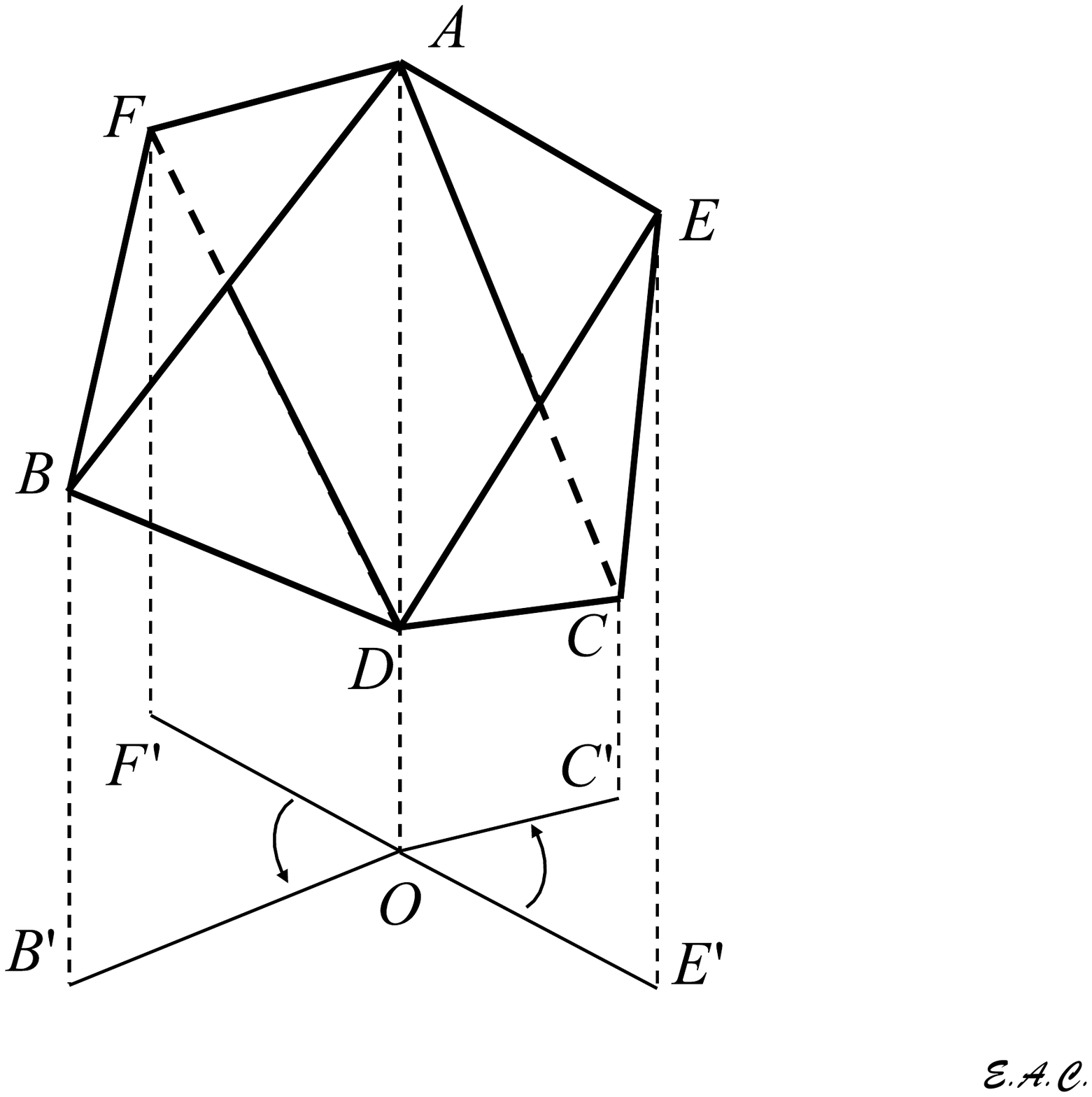,width=6.0in}}
\caption{}
\label{fig:flex1}
\end{figure}
The figure formed by the two last triangles is evidently superimposable 
either onto the figure formed by the first two, or to the one symmetric
to this with respect to an arbitrary plane. This leads to two cases that
must be examined.

Let us suppose then that the two systems of triangles form superimposable
figures. Now the four triangles $\triangle ADF$, $\triangle ADB$, 
$\triangle ADE$, $\triangle ADC$ are projected,
on a plane perpendicular to $AD$, respectively along the lines
$OF'$, $OB'$, $OE'$, $OC'$; such that the angles $\angle F'OB'$, 
$\angle F'OC'$ will be
equal and that the same sense of rotation brings $OF'$ inot coincidence with
$OB'$, $OE'$ into coincidence with $OC'$. One then has
\[ \angle F'OE' = \angle B'OC'\ . \]
From this there results that the two trihedrals $A(FDE)$ and $D(CAB)$,
that have the same orientation, are equal by virtue of having an equal
dihedral comprised of two faces equal each to each. In effect, the preceding
equality expresses that property of those dihedrals having $AD$ as a common
edge. One has, on the other hand,
\[ \angle FAD = \angle CDA\ , \]
\[ \angle DAE = \angle ADB\ . \]
One is led from this to the equality of the angles $\angle FAE$, $\angle CDB$. 
The two triangles $\triangle FAE$, $\triangle BDC$ are therefore equal, 
and one has
\[ FE = BC\ .\]
This equality is true during all the deformations of which our system 
of triangles is susceptible, under the condition, I repeat, that the set of
the latter two is at all times superimposable to that of the former two.

Now, this deformation is such that the complete determination of the system
depends on two parameters (for which one may take for example the distance
$AD$ and the angle $\angle F'OE'$). This deformation is thus still possible if the
system is subjected to the supplementary condition that the distance $EF$
remains constant: it will then be the same as the distance $BC$.

The set of the figure will exhibit eight invariant triangles, constituting
a defomable octahedron whose opposing edges are pairwise equal.

This octahedron admits one axis of symmetry: let us draw, in effect, in the 
bisecting plane of the dihedral projected on $B'OE'$, a straight line $L$,
perpendicular to $AD$ and passing through the midpoint of that segment.
It is clear that the points $A$, $B$, $F$ are respectively symmetric to
the points $D$, $E$, $C$, with respect to $L$. On may then bring the
octahedron into coincidence with itself by making it turn by two right
angles about $L$. On may see also that the three diagonals of the octahedron
are perpendicular to the same line $L$ which divides each of them in two
equal parts \cite{Mannheim}.

One may construct such an octahedron by means of cardboard
triangles appropriately cut and assembled with sticky tape.
It is necessary to leave empty
the facets $ABC$, $DEF$, which are only realised by their contour.

\begin{figure}
\centerline{\epsfig{figure=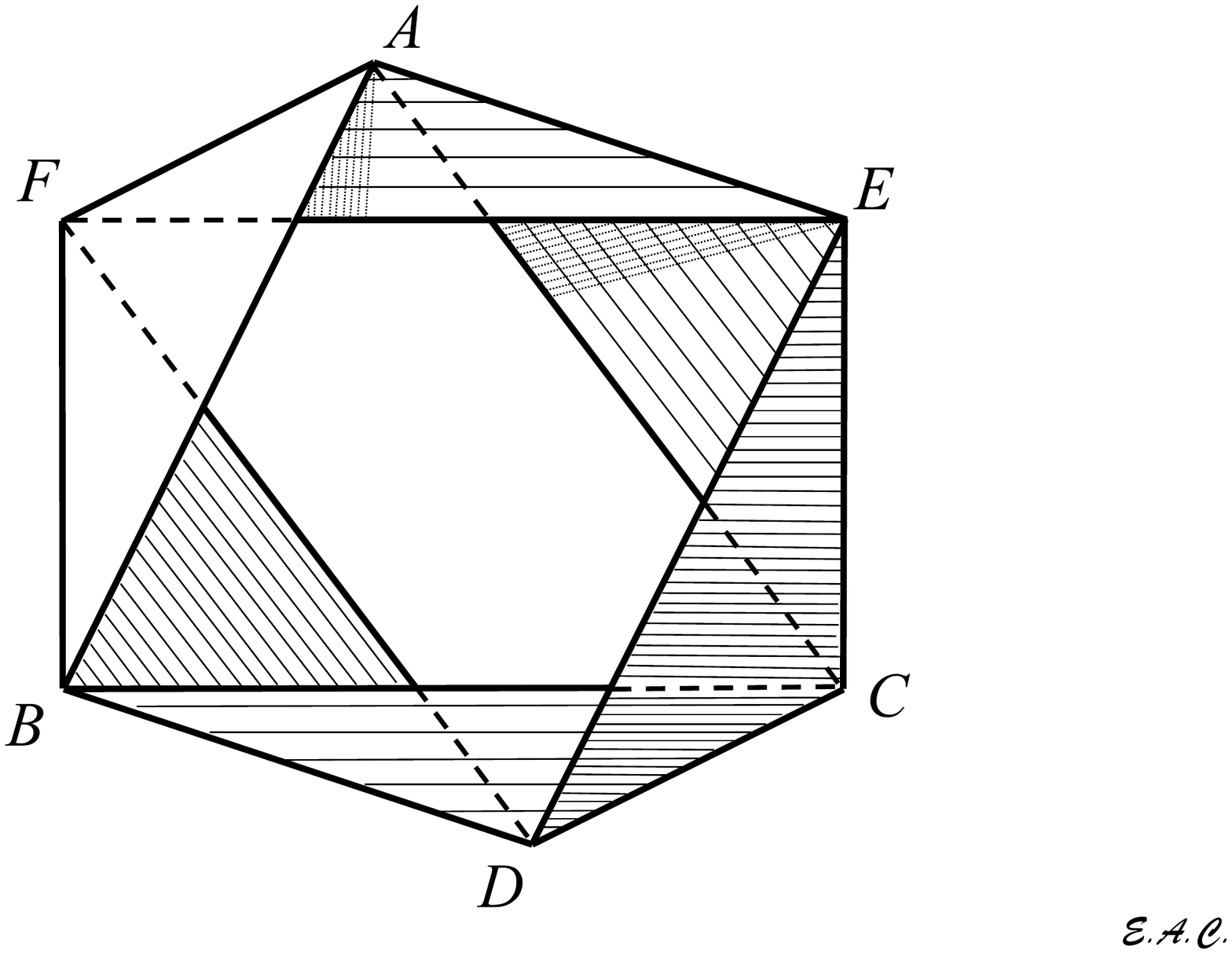,width=6.0in}}
\caption{Octahedron whose opposite edges are equal in pairs. The facets are
$ABC$, $DEF$, $BCD$, $CAE$, $ABF$, $AEF$, $BFD$, $CDE$}
\label{fig:octahedron1}
\end{figure}

The model thus obtained is represented in Fig. \ref{fig:octahedron1}.

If, now, returning to Fig. \ref{fig:flex1}, one assumes that the figure 
formed by the two triangles $\triangle ACE$, $\triangle DCE$, 
is superimposable to the figure
symmetric to that formed by the two triangles $\triangle AFB$, $\triangle DFB$, 
one will
conclude, by an argument entirely similar to the preceding one that, if the 
distance $EF$ remains constant, the distance $BC$ is necessarily variable:
the octahedron $ABCDEF$ is not deformable. This establishes the assertion
made at the end of section VI. One sees that the relations 
(\ref{typeI}) are not independent, but that one of them implies the other 
five. 
\section{(The case of a rhomboidal tetrahedral)}
I pass to the case where one among the tetrahedral angles is rhomboidal,
so that only one of them has its adjacent faces equal or supplementary
pairwise. Let us suppose that the tetrahedral angle $C$ has this
special property, that is it has, for example:
\[\angle BCD = \angle DCE\,,\qquad
  \angle BCA = \angle ACE\,. \]
The system (\ref{tetra}), (\ref{tetrb}), (\ref{tetrc}), becomes then
\begin{equation}
\label{tetra1}
At^2u + 2Ct +Du = 0 \ ,
\end{equation}
\begin{equation}
\label{tetrb1}
A^\prime t^2v^2 + B^\prime t^2 +2C^\prime tv +D^\prime v^2 + E^\prime  = 0 \ ,
\end{equation}
\begin{equation}
\label{tetrc1}
A^{\prime\prime} u^2v^2 + B^{\prime\prime} u^2 +2C^{\prime\prime} uv 
+D^{\prime\prime} v^2 + E^{\prime\prime}  = 0 \ .
\end{equation}
We may reason as before: the equations (\ref{tetrb1}) and (\ref{tetrc1})
can have at all times the two roots in $v$ common, or instead they may 
only have one of these roots in common.

I will demonstrate shortly that the first assumption is inadmissible.
The equations (\ref{tetrb1}) and (\ref{tetrc1}) have thus a single
common root in $v$; this root is a rational function of $t$ and $u$,
and in terms of $t$, equation (\ref{tetrb1}) reduces necessarily
to one of the forms
\[ A^\prime t^2v + 2C^\prime t +D^\prime v  = 0 \ ,\ 
   B^\prime t^2 +2C^\prime tv + E^\prime  = 0 \ .
\]
In other words, the tetrahedral angle $B$ is rhomboidal, and it has
\[ \mbox{dihedral}\ BF = \mbox{dihedral}\ BC = \mbox{dihedral}\ CE \ .\]
It also follows that the equations (\ref{tetra1}) and (\ref{tetrc1}),
respectively in $t$ and in $v$, have double roots for the same values
of $u$. One then concludes that the dihedrals $CD$ and $AF$ remain equal or
supplementary at all times. The same is true for the dihedrals $BD$, $AE$.
Continuing in this way, it will be true that the opposing pairs of the 
dihedrals of the octahedron all possess the same property.\\
\indent The conclusions of Section VI are not therefore changed.\\
\indent It remains to show that the equations (\ref{tetrb1}) and (\ref{tetrc1})
cannot have their two roots in $v$ identically equal at all times.\\
\indent In effect, if it is so, one must have
\[
\frac{A't^2+D'}{A''u^2+D''} =\frac{C't}{C''u} =
\frac{B't^2+E'}{B''u^2+E''}\]
or
\begin{eqnarray*}
 A'C''t^2u-C'A''tu^2-C'D''t+D'C''u & = & 0\ ,\\
 B'C''t^2u-C'B''tu^2-C'E''t+E'C''u & = & 0\ .
\end{eqnarray*}
These last equations must be identical with equation (\ref{tetra1}).
One then has
\[ C'A'' = 0\ ,\ C'B'' = 0 \ ;\]
from which
\[ C' = 0 \ \ \mbox{or instead}\ \ A''= 0\ ,\ \ B'' = 0 \ . \]
Now, if it is true that
\[ C' = 0\ ,\]
there follows from the discussion of Sec. II that equation (\ref{tetrb1})
has of necessity one of the forms
\[ A't^2v^2 + E' = 0\ ,\ B't^2+D'v^2 = 0\ .\]
The tetrahedral angle $B$ will then have its opposing faces equal or 
supplementary pairwise, and we have excluded such tetrahedral angles 
in examining the present case.\\
\indent One then has
\[ A''= 0\ ,\ \ B'' = 0 \ , \]
and equation (\ref{tetrc1}) reduces to
\[ 2C''uv + D'' v^2 + E'' = 0 \ .\]
The tetrahedral angle $A$ is rhomboidal (Sec.II), and one has
\[ \angle EAC = \pi - \angle BAC\ ,\ \ 
   \angle FAE = \pi - \angle BAF\ ,\]
\[ \mbox{dihedral} AE = \mbox{dihedral} AB\ .\]
Since one also has
\[ \mbox{dihedral} CE = \mbox{dihedral} CB\ ,\]
the two tetrahedral angles having their vertices at $B$ and $E$ must
deform in such a way that two dihedrals adjacent to the first are
at all times equal, respectively, to two dihedrals adjacent to the 
second tetrahedral angle. It must also be true (Sec. III) that
their faces are pairwise equal or supplementary. One has also
\[ \mbox{dihedral} FE = \mbox{dihedral} FB\ ,\]
\[ \mbox{dihedral} DE = \mbox{dihedral} DB\ .\]
The tetrahedral angles $F$ and $D$ are also rhomboidal. Combining these
results, one can see that the triangles
\[ \triangle AFE\ \ \ \mbox{and}\ \ \ \triangle AFB\ ,\]
\[ \triangle ACB\ \ \ \mbox{and}\ \ \ \triangle ACE\ ,\]
\[ \triangle BFD\ \ \ \mbox{and}\ \ \ \triangle EFD\ ,\]
\[ \triangle BCD\ \ \ \mbox{and}\ \ \ \triangle ECD\  \]
have, pairwise, their angles equal or supplementary. The equality alone
is possible. One has in particular
\[ \angle BAF = \angle EAF\ ,\]
\[ \angle BAC = \angle EAC\ .\]
We have concluded, on the other hand, that these angles, faces of the
tetrahedral angle $A$, are pairwise supplementary. The may not all be also 
equal unless they are all $\pi/2$, which is visibly impossible.

\section{(The case of a unicursal tetrahedral)}
There remains to examine the case where at least one of the tetrahedral 
angles has its opposing faces equal or supplementary in pairs.\\
\indent Let us suppose this is so for the tetrahedral angle $C$. The
variables $t$ and $u$ satisfy one of the relations
\[ tu = k\ ,\ \frac{t}{u}=k\ . \]
I shall admit the existence of the first, the reasoning being the same for the 
second case. Assuming this, two distinct hypotheses need to be examined:\\
\begin{enumerate}
\item None of the tetrahedral angles whose vertices are located at the 
ends of the edges issuing from $C$ has its opposing faces equal or 
supplementary.
\item One of these tetrahedral angles prossesses that property.
\end{enumerate}

In the first case, the tetrahedral angles $A$ and $B$ are general or rhomboidal.
Let us suppose them general, for example. The system of relations among 
$t, u, v$ is then
\begin{equation}
\label{tetra2}
tu = k \ ,
\end{equation}
\begin{equation}
\label{tetrb2}
A^\prime t^2v^2 + B^\prime t^2 +2C^\prime tv +D^\prime v^2 + E^\prime  = 0 \ ,
\end{equation}
\begin{equation}
\label{tetrc2}
A^{\prime\prime} u^2v^2 + B^{\prime\prime} u^2 +2C^{\prime\prime} uv 
+D^{\prime\prime} v^2 + E^{\prime\prime}  = 0 \ .
\end{equation}
Equation (\ref{tetrb2}) in $t$ and equation (\ref{tetrc2}) in $u$ have their 
roots equal for the same values of $v$: it then follows that the dihedrals
$AE$ and $BD$ are at all times equal or supplementary.\\
\indent Concurrent consideration of the tetrahedral angles $A, C, E$ shows
similarly that the dihedrals $AB$, $DE$ remain at all times equal or 
supplementary. One then may construct fig. \ref{fig:triangle2}, where the
dihedrals marked by the same number exhibit the same relationship (the dihedrals
$1$ and the dihedrals $2$ due to the nature of the tetrahedral angle $C$).\\
\indent If one applies to the tetrahedral angles $A$ and $D$ on the one hand,
$B$ and $E$ on the other, the theorem which has already been utilised 
several times, one sees:
\begin{enumerate}
\item That the dihedrals $FA$ and $FD$, on the one hand, $FB$ and $FE$ on the 
other, are at all times equal, and that the tetrahedral angle $F$ has,
consequently, its opposing faces equal or supplementary in pairs;
\begin{figure}
\centerline{\epsfig{figure=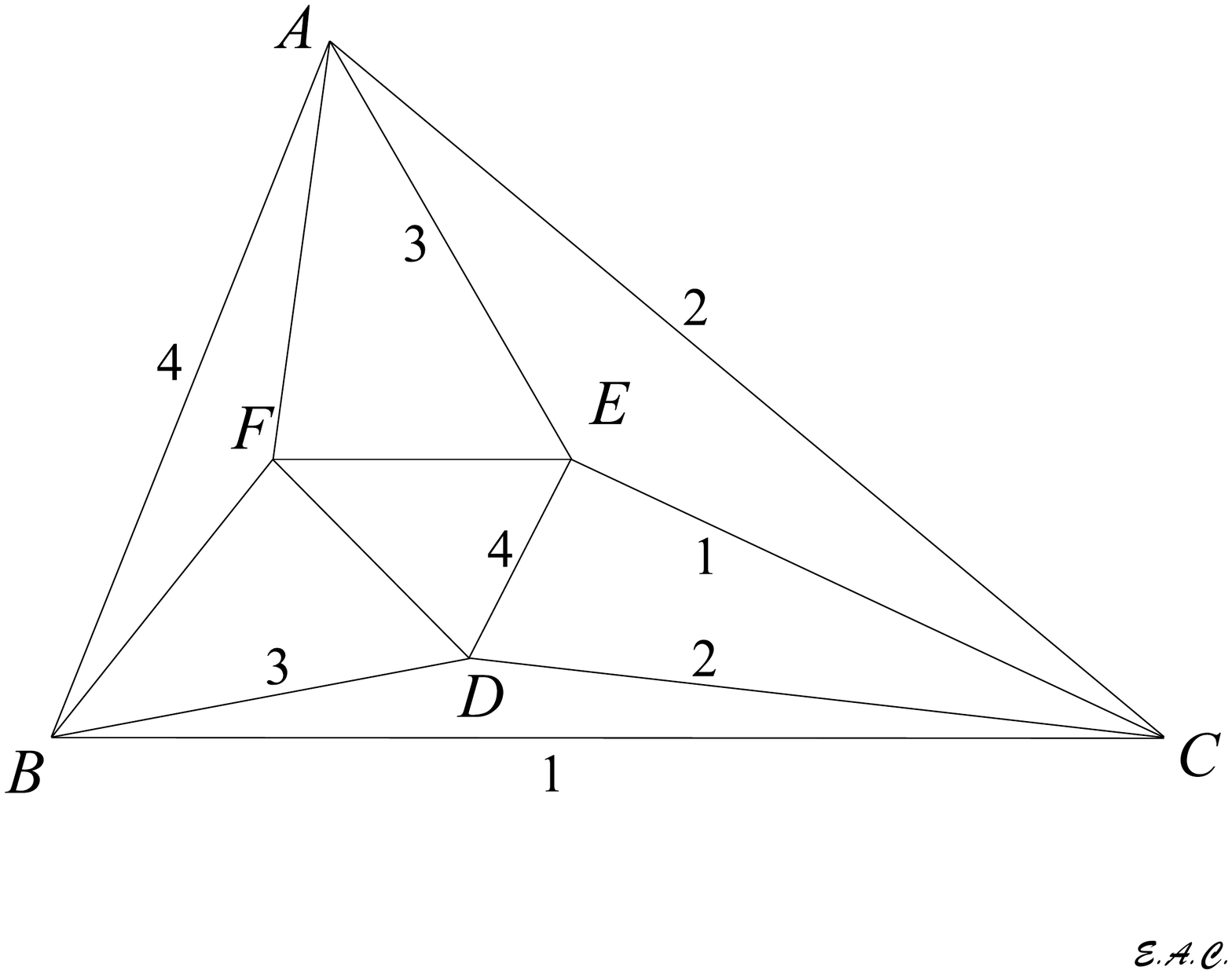,width=4.4in}}
\caption{}
\label{fig:triangle2}
\end{figure}
\item That the facets of the octahedron are subdivided into four pairs of 
triangles having their angles equal in pairs:
\end{enumerate}
\[ \triangle AFE \ \ \mbox{and}\ \ \triangle DFB, \]
\[ \triangle AFB \ \ \mbox{and}\ \ \triangle DFE, \]
\[ \triangle AEC \ \ \mbox{and}\ \ \triangle DBC, \]
\[ \triangle ABC \ \ \mbox{and}\ \ \triangle DEC. \]
For each pair the homologous vertices are written in the same order. One has
therefore the sequence of equalities
\[ \frac{AF}{DF} = \frac{FE}{FD} = \frac{EA}{BD}\ ,\ \ 
   \frac{AF}{DF} = \frac{FB}{FE} = \frac{BA}{ED}\ ,\]
\[ \frac{AE}{DB} = \frac{EC}{BC} = \frac{BA}{ED}\ ,\ \ 
   \frac{AB}{DE} = \frac{BC}{EC} = \frac{CA}{CD}\ ,\]
from which it is seen that
\begin{equation}
\label{typeII}
\left\{
\begin{array}{rrr}
FA = FD\ , & FE = FB\ , & CA = CD\ ,\\
CB = CE\ , & AE = DB\ , & AB = DE\ .
\end{array}
\right.
\end{equation}
{\em Therefore the edges of the octahedron must again be pairwise equal, but 
the equal edges are not all pairwise opposite, a fact that distinguishes 
conditions (\ref{typeII}) from conditions (\ref{typeI}).}

\section{(Flexibility of Type II: Opposing unicursal tetrahedrals
and a plane of symmetry)}
\begin{figure}
\centerline{\epsfig{figure=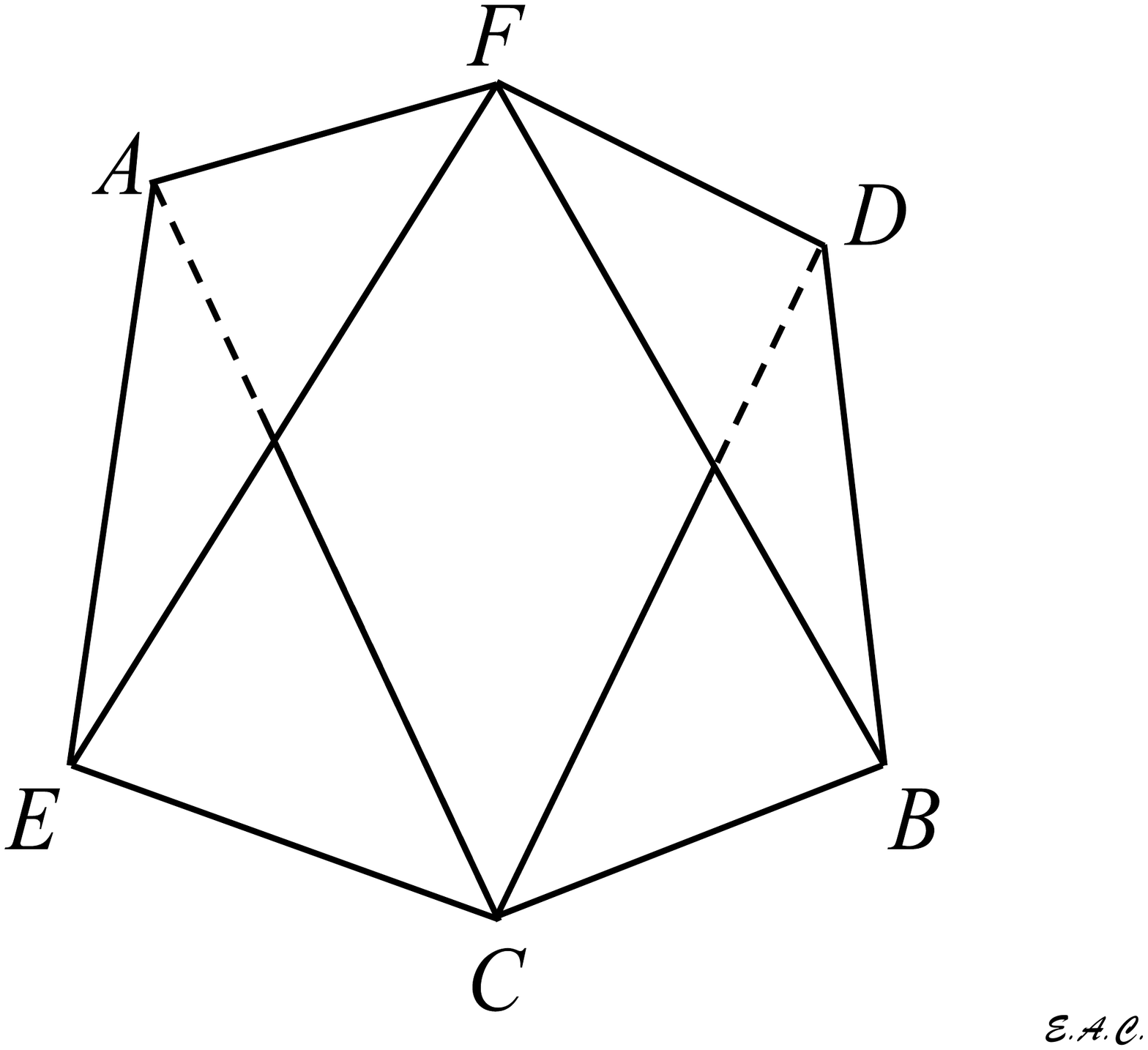,width=4.2in}}
\caption{}
\label{fig:flex2}
\end{figure}

\begin{figure}
\centerline{\epsfig{figure=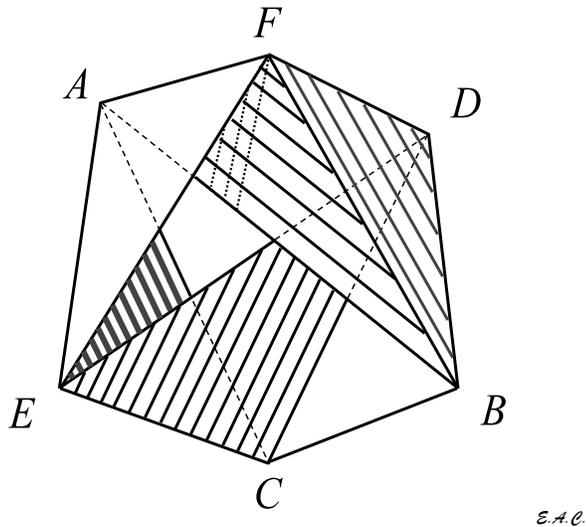,width=4.1in}}
\caption{Octahedron possessing a plane of symmetry passing through two
opposite vertices. The facets are $ABC$, $DEF$, $BCD$, $CAE$, $ABF$,
$AEF$, $BFD$, $CDE$.}
\label{fig:octahedron2}
\end{figure}

It must now be shown that all octahedra whose edges satisfy relations
(\ref{typeII}) are deformable, provided certain conditions, related to the
placement of the facets, are satisfied.\\
\indent For this, let us consider (fig. \ref{fig:flex2}) the system of four
rigid triangles
$\triangle FAE$, $\triangle CAE$, $\triangle FDB$, $\triangle CDB$, 
joined at the points $F$, $C$, and along the
lines $AE$, $BD$. Assume the equalities
\[ FA=FD, FE=FB, CA=CD, CE=CB, AE=DB.\]
The system of the last two triangles is superimposample on that of the first
two, or instead it is symmetric to that with respect to a plane passing through
$FC$. \\
\indent In the first case it will be true, by arguments analogous to those in
Sec. VII, that one cannot have at all times that $AB=DE$. The octahedron
formed by connecting the points $A$ and $B$ on the one hand, $D$ and $E$ on the
other, not satisfying all of the relationships (\ref{typeII}) and 
moreover not satisfying relations (\ref{typeI}), is not deformable.\\
\indent On the contrary, if the systems of triangles are symmetric, it is 
evident that, whatever their relative position, one has always
\[ AB = DE\ .\]
The reasoning will proceed like that in Sec. VII, and it will establish
that relations (\ref{typeII}), just like relations (\ref{typeI}),
suffice, under the indicated restrictions, to ensure the deformability of an
octahedron.\\
\indent This second octahedron may be realised like the first, by leaving
empty the facets $ABC$ and $DEF$. The model thus obtained is represented in 
fig. \ref{fig:octahedron2}.

\section{Flexibility of Type III: adjacent unicursal tetrahedrals in proportion)}
I arrive finally at the case where two tetrahedral angles having their 
vertices adjacent, the tetrahedral angle $C$ and the tetrahedral angle
$B$, for example, have each of their opposite faces equal or supplementary
in pairs.\\
\indent One can see immediately that this must be similarly true for all
the tetrahedral angles of the octahedron.\\
\indent In effect, from the relations
\[ tu \ \ \mbox{or} \ \ \frac{t}{u} = k,\]
\[ tv \ \ \mbox{or} \ \ \frac{t}{v} = k^\prime,\]
one shows 
\[ \frac{u}{v} \ \ \mbox{or} \ \ uv = k^{\prime\prime},\]
which establishes the proposition for the tetrahedral angle $A$; this will be
established similarly for the other tetrahedral angles $D$, $E$, $F$.\\
\indent If then the dihedral $BD$ becomes equal to $0$ or $\pi$, $t$ becoming
zero or infinity, the variables $u$, $v$ are also zero or infinite, and the
dihedrals $AC$, $AB$, equal to $0$ or $\pi$. It is similarly so for all the 
other dihedrals. In other words, the octahedron may be completely flatened
onto the facet $ABC$.\\
\indent Let us represent it in that position. There may be several cases of this
figure, for which the reasoning is identical.

\begin{figure}
\centerline{\epsfig{figure=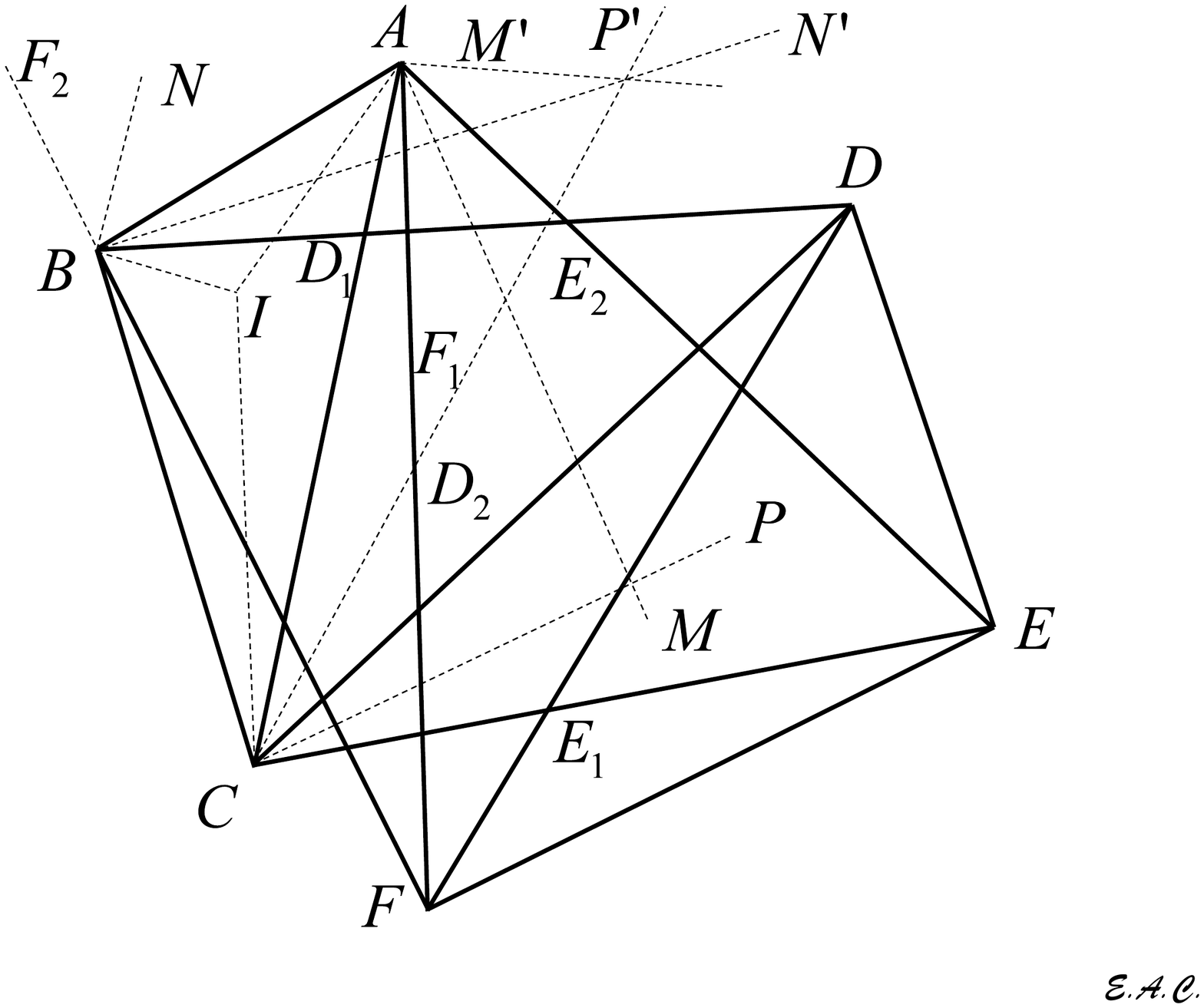,width=6.0in}}
\caption{}
\label{fig:flex3a}
\end{figure}
I will suppose, for example, that the arrangement is that 
of fig. \ref{fig:flex3a}. One has
\begin{eqnarray*}
\angle FAE & = & \angle BAC\ ,\\
\angle DCE & = & \angle ACB\ ,\\
\angle DBF & = & \pi - \angle ABC\ .\\
\end{eqnarray*}
During the deformation of the octahedron, one may have various systems
of relationships among $t,u,v$. For example (sec. II), let
\begin{eqnarray*} 
uv & = & \sfrac{\cos{\sfrac{\angle BAF-\angle BAC}{2}}}
{\cos{\sfrac{\angle BAF+\angle BAC}{2}}}\ ,\\
\frac{t}{v} & = & \sfrac{\sin{\sfrac{\angle ABC-\angle DBC}{2}}}
{\sin{\sfrac{\angle ABC+\angle DBC}{2}}}\ ,\\
tu & = & \sfrac{\sin{\sfrac{\angle DCB-\angle ACB}{2}}}
{\sin{\sfrac{\angle DCB+\angle ACB}{2}}}\ .
\end{eqnarray*}
These equations must be satisfied by an infinity of sets of values of $t,u,v$.
One has then
\begin{equation}
\label{typeIII}
\sfrac{\cos{\sfrac{\angle BAF-\angle BAC}{2}}}
{\cos{\sfrac{\angle BAF+\angle BAC}{2}}}
\sfrac{\sin{\sfrac{\angle ABC-\angle DBC}{2}}}
{\sin{\sfrac{\angle ABC+\angle DBC}{2}}}
\sfrac{\sin{\sfrac{\angle DCB+\angle ACB}{2}}}
{\sin{\sfrac{\angle DCB-\angle ACB}{2}}} = 1\ .
\end{equation}
This is the necessary and sufficient condition for the octahedron
$ABCDEF$ to be deformable.\\
\indent In order to give to this condition a geometric form, let us draw the
{\em interior} bisectrices of the triangle 
$\triangle ABC$: $AI$, $BI$, $CI$. Let us draw 
also the lines $AM$, $CP$, $BN$, the first two, {\em interior} bisectrices
of the angles $\angle EAF$, $\angle DCE$, the third, {\em exterior} 
bisectrice of the angle $\angle FBD$. Let finally $FM^\prime$ be 
the {\em exterior} bisectrice of the angle $\angle IAM$, $BN^\prime$ 
and $CP^\prime$ the {\em exterior} bisectrices
to the angles $\angle IBN$, $\angle ICP$. One has:
\begin{eqnarray*}
\sfrac{\cos{\sfrac{\angle BAF-\angle BAC}{2}}}{\cos{\sfrac{\angle BAF+\angle BAC}{2}}} & = &
\sfrac{\cos\left(\sfrac{\angle IAM}{2}-\angle IAC\right)}{\cos\left(\sfrac{\angle IAM}{2}+\angle BAI\right)} =
\sfrac{\cos{\left(\angle IAM^\prime-\sfrac{\pi}{2}-\angle IAC\right)}}
    {\cos{\left(\angle IAM^\prime-\sfrac{\pi}{2}+\angle BAI\right)}} \\
&=& \sfrac{\sin{\left(\angle IAM^\prime - \angle IAC\right)}}
          {\sin{\left(\angle IAM^\prime + \angle BAI\right)}} 
   = \sfrac{\sin{\angle CAM^\prime}}{\sin{\angle BAM^\prime}}\ ,
\end{eqnarray*} 
\[
\sfrac{\sin{\sfrac{\angle ABC-\angle DBC}{2}}}{\sin{\sfrac{\angle ABC+\angle DBC}{2}}} = 
\sfrac{\sin{\left(\angle IBA - \sfrac{\angle IBN}{2}\right)}}
      {\sin{\left(\angle CBI + \sfrac{\angle IBN}{2}\right)}} =
\sfrac{\sin{\left(\angle IBA - \angle IBN^\prime\right)}}
      {\sin{\left(\angle CBI - \angle IBN^\prime\right)}}     =
\sfrac{\sin{\angle N^\prime BA}}{\sin{\angle CBN^\prime}}\,
\]
\[
\sfrac{\sin{\sfrac{\angle DCB+\angle ACB}{2}}}{\sin{\sfrac{\angle DCB-\angle ACB}{2}}} =
\sfrac{\sin{\left(\sfrac{\angle PCI}{2} + \angle ICB\right)}}
      {\sin{\left(\sfrac{\angle PCI}{2} - \angle ACI\right)}} =
\sfrac{\sin{\left(\angle P^\prime CI + \angle ICB\right)}}
      {\sin{\left(\angle P^\prime CI - \angle ACI\right)}}     =
\sfrac{\sin{\angle P^\prime CB}}{\sin{\angle P^\prime CA}}\ .
\]
Relation (\ref{typeIII}) becomes then
\[ 
\sfrac{\sin{\angle CAM^\prime}}{\sin{\angle BAM^\prime}}
\sfrac{\sin{\angle N^\prime BA}}{\sin{\angle CBN^\prime}}
\sfrac{\sin{\angle P^\prime CN}}{\sin{\angle P^\prime CA}} = 1\ ; \]
from which it follows, by virtue of a well known theorem, that the lines
$AM^\prime$,$BN^\prime$,$CP^\prime$, meet at a point (are concurrent).\\
\indent I have made a particular hypothesis on the form of the relationships
that exist among $t,u,v$. It is clear that in all of the other cases, one 
will arrive at a similar result, and it may be pronounced as the general rule
for constructing an artriculated octahedron all of whose tetrahedral 
angles have their opposite faces pairwise equal or supplementary.

{\em Construct an arbitrary triangle $\triangle ABC$, 
whose interior bisectrices are
the lines $AI$, $BI$, $CI$, and from the vertices of that triangle draw three
concurrent lines $AM^\prime$, $BN^\prime$, $CP^\prime$. Trace the lines
$AM$, $BN$, $CP$, respectively symmetric to the lines $AI$, $BI$, $CI$, 
with respect to the lines $AM^\prime$, $BN^\prime$, $CP^\prime$. \\
\indent Construct then the angles $\angle F_1AE_2$, $\angle D_1BF_2$, 
$\angle E_1CD_2$, obtained by
rotating in their own planes the angles 
$\angle BAC$, $\angle CBA$, $\angle ACB$, about their 
vertices by angles equal in magnitude and sign respectively to 
$\angle IAM$, $\angle IBN$, $\angle ICP$. 
Let $D$, $E$, $F$ be the points of intersection, respectively, of the
lines (suitably elongated through the vertices of the triangle) 
$BD_1$ and $CD_2$, $CE_1$ and $AE_2$, $AF_1$ and $BF_2$.\\
\indent Assuming the triangles $\triangle ABC$, $\triangle BCD$, 
$\triangle CAE$, $\triangle ABF$, 
$\triangle AEF$, $\triangle BFD$, $\triangle CDE$, $\triangle DEF$, 
are constructed and articulated in pairs along their common edges, these
triangles are the facets of a deformable octahedron.}

\begin{figure}
\centerline{\epsfig{figure=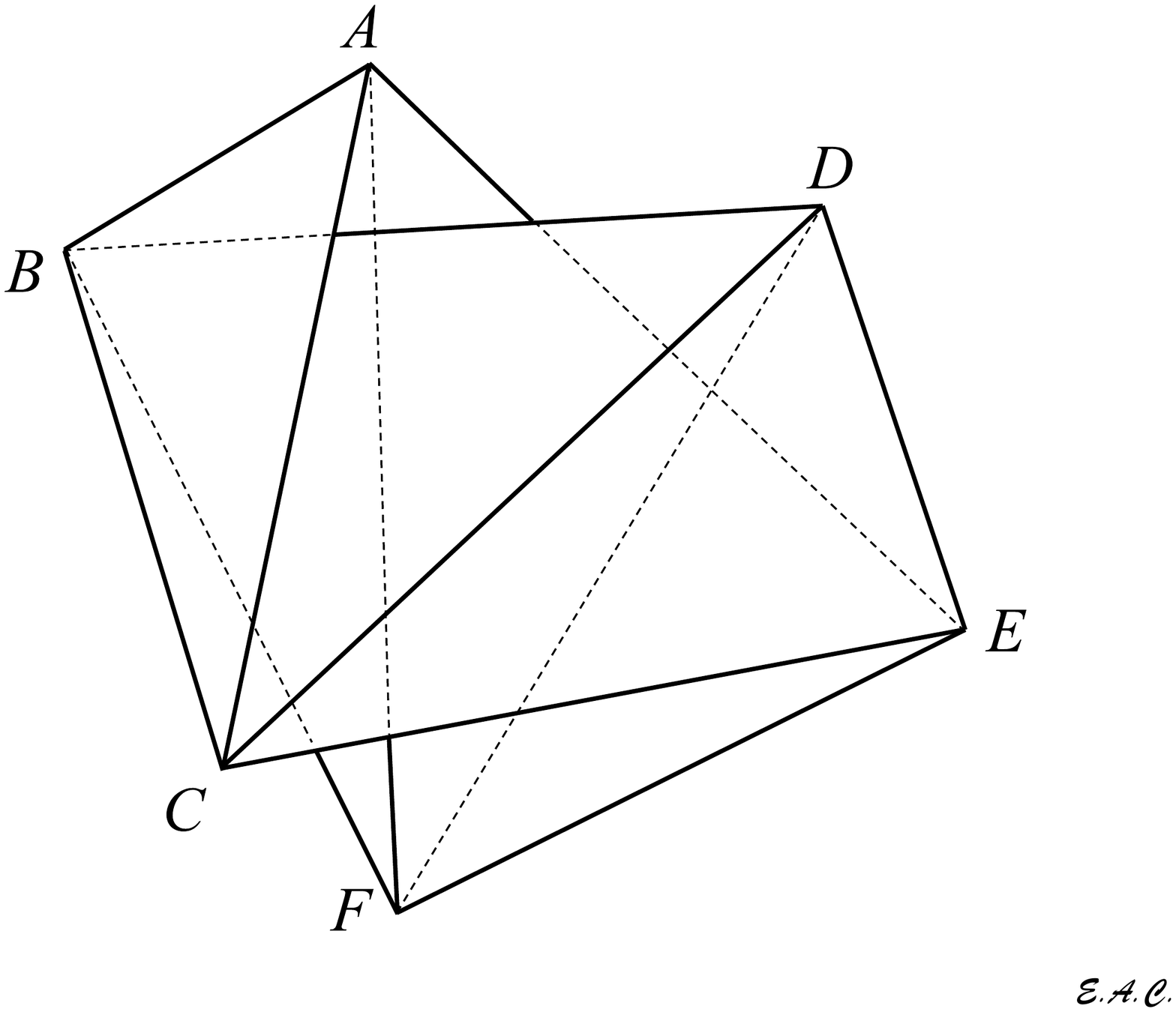,width=4.2in}}
\caption{}
\label{fig:flex3b}
\end{figure}

Certain difficulties may be encountered during the construction of a model
of this latter octahedron: which facets should be left open? from which
side should the concavity of each dihedral be oriented? These are problems
that may not be resolved in each particular case without a careful examination
of the relations among $t, u, v$ and the signs imposed on each variable.
It would be too lengthy to discuss here the detailed reasoning. I shall content
myself by indicating how to carry out the construction of the octahedron
of fig. \ref{fig:flex3a}.\\
\indent I have represented it (fig. \ref{fig:flex3b}, right) in the same 
position
as in the preceding figure; the facets $AEC$, $DBF$ are empty. The dihedral
$AF$ has its concavity facing forward from the plane of the figure; that of
the dihedral $DC$ is facing backward.


\begin{figure}
\centerline{\epsfig{figure=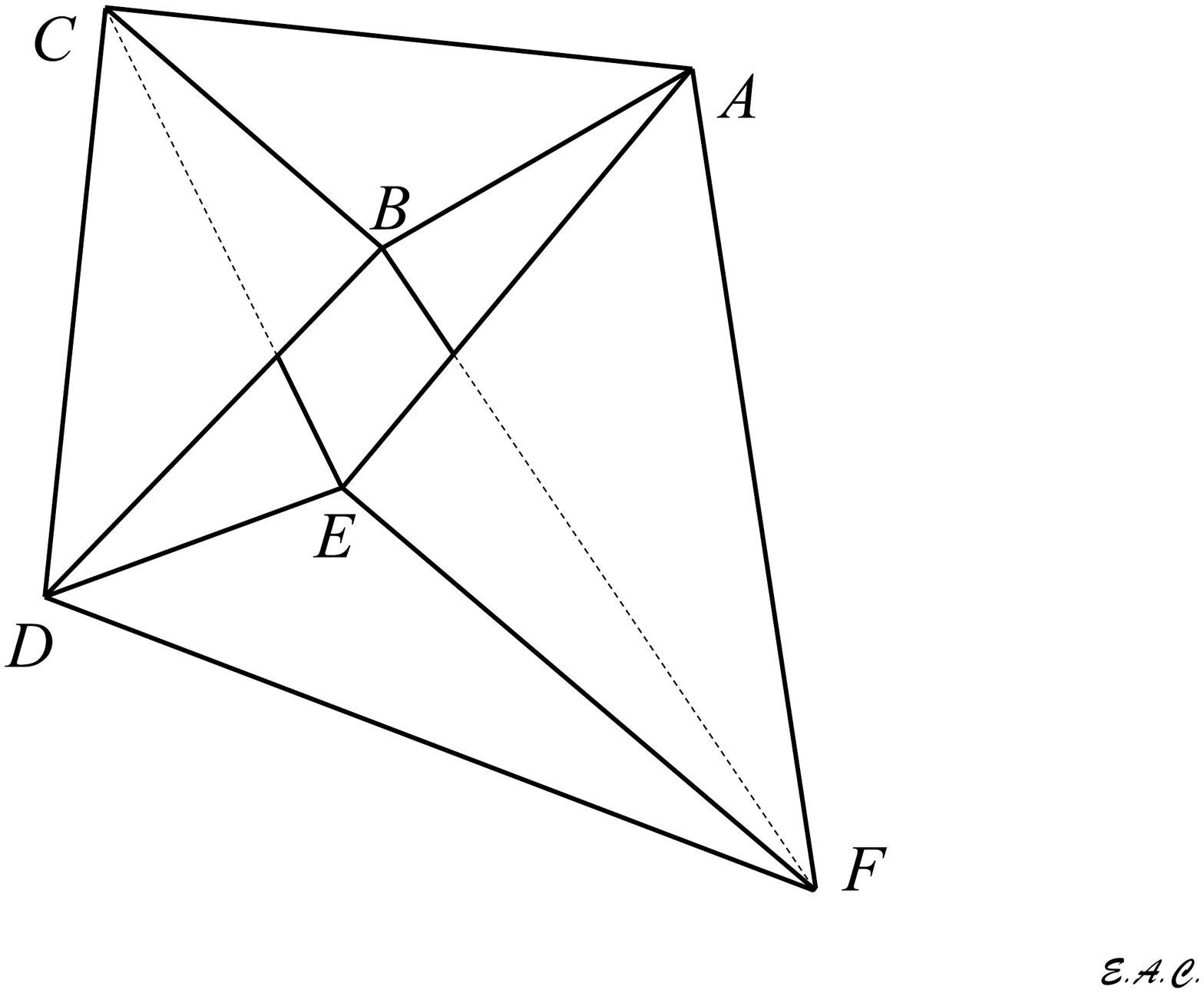,width=4.2in}}
\caption{}
\label{fig:flex3c}
\end{figure}
\begin{figure}
\centerline{\epsfig{figure=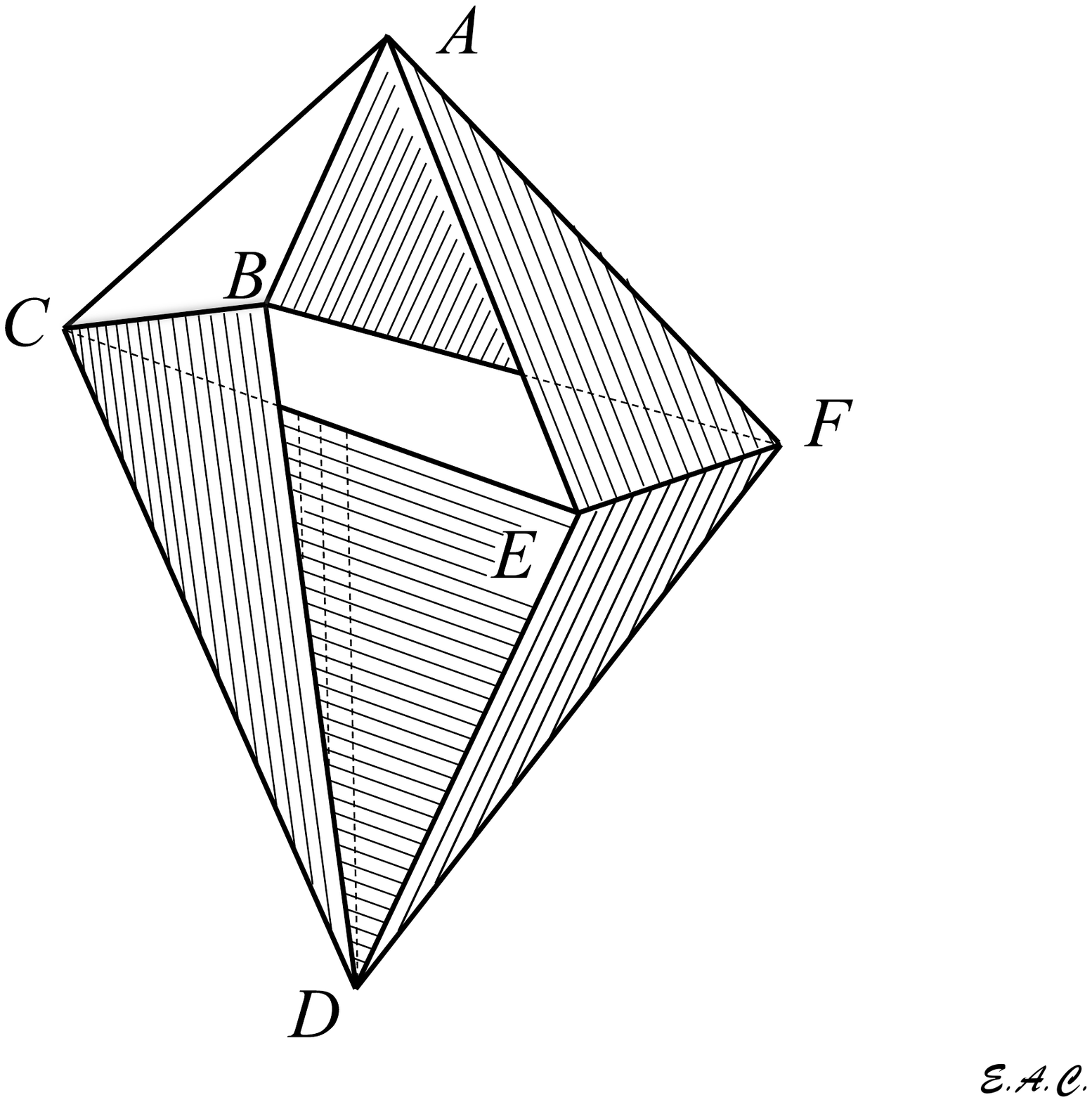,width=4.4in}}
\caption{Octahedron all of whose tetrahedral angles have their opposing faces
pairwise equal or supplementary.
The facets are $ABC$, $DEF$, $BCD$, $CAE$, $ABF$, $AEF$, $BFD$, $CDE$.}
\label{fig:octahedron3}
\end{figure}

The broken lines indicate unambiguously in which order they are
superimposed on the facets: in particular, the facet $DEF$ is behind the
facet $AEF$. \\
\indent When deforming this octahedron, by setting the facet $ABF$ on the
plane of the figure, the vertices $D$, $E$, $C$ are displaced in front of this
plane, and pursuing the deformation one arrives at a new flattened position 
(fig. \ref{fig:flex3c}). Fig. \ref{fig:octahedron3} represents an
intermediate position. \\
\indent If the reader takes the pains of constructing this model, 
by following the preceding instructions, I shall reiterate that the proportions
of the figure must be preserved in the most precise fashion.

\section{(Discussion)}
In summary, the preceding study has shown that there exist {\it three} types
of articulated octahedra with invariant facets. All these polyhedra are 
concave or, to be more precise, possess certain facets that intercross.

Octahedra of type I (\ref{typeI}) and II (\ref{typeII}) have simple 
definitions: the first possess an 
axis of symmetry and, as a result, are such that the figure formed by four
of their facets with a common vertex is superimposable on the figure formed
by the other four; those of type II have a plane of symmetry, passing through
two opposite vertices. (These definitions are not {\it absolutely}
sufficient, but one can only complete them by lengthy discussion, and
I feel that an examination of figures \ref{fig:octahedron1} 
and \ref{fig:octahedron2} makes that unnecessary).

About octahedra of type III (\ref{typeIII}), 
we have seen that their definition is more
complicated; their deformability is far from being as intuitive as that of
the former types, and in this sense they ought to be considered as the most
interesting.

I shall also remark that the problem which I have discussed is identical to the
following two problems:

$1^o$ {\it What are the deformable oblique hexagons with constant sides and
angles?}

If, in fact, an oblique hexagon is deformable under these conditions,
the segments that join its next to nearest vertices have constant lengths and 
the eight triangles formed by these segments and by the sides of the
hexagon are the facets of a deformable octahedron with constant edges.

A deformable octahedron exhibits, on the other hand, four such hexagons.
These are (fig. \ref{fig:octahedron1}, \ref{fig:octahedron2} or 
\ref{fig:flex3c}) the hexagons
\[ ABCDEF\ ,\]
\[ ABFDEC\ ,\]
\[ AECDBF\ ,\]
\[ AEFDBC\ .\]
$2^o$ {\it Under what conditions is a system of six planes 1, 2, 3, 4, 5, 6,
where each is articulated with the previous along a line serving as a hinge,
with plane 6 being articulated with planes 1 and 5, susceptible to
deformation?}

In fact, the lines along which these planes are articulated form a deformable
hexagon with fixed sides and angles.

The planar facets of one of the octahedra of figs. \ref{fig:octahedron1}, 
\ref{fig:octahedron2}, \ref{fig:octahedron3} constitute such a system.
\begin{figure}
\centerline{\epsfig{figure=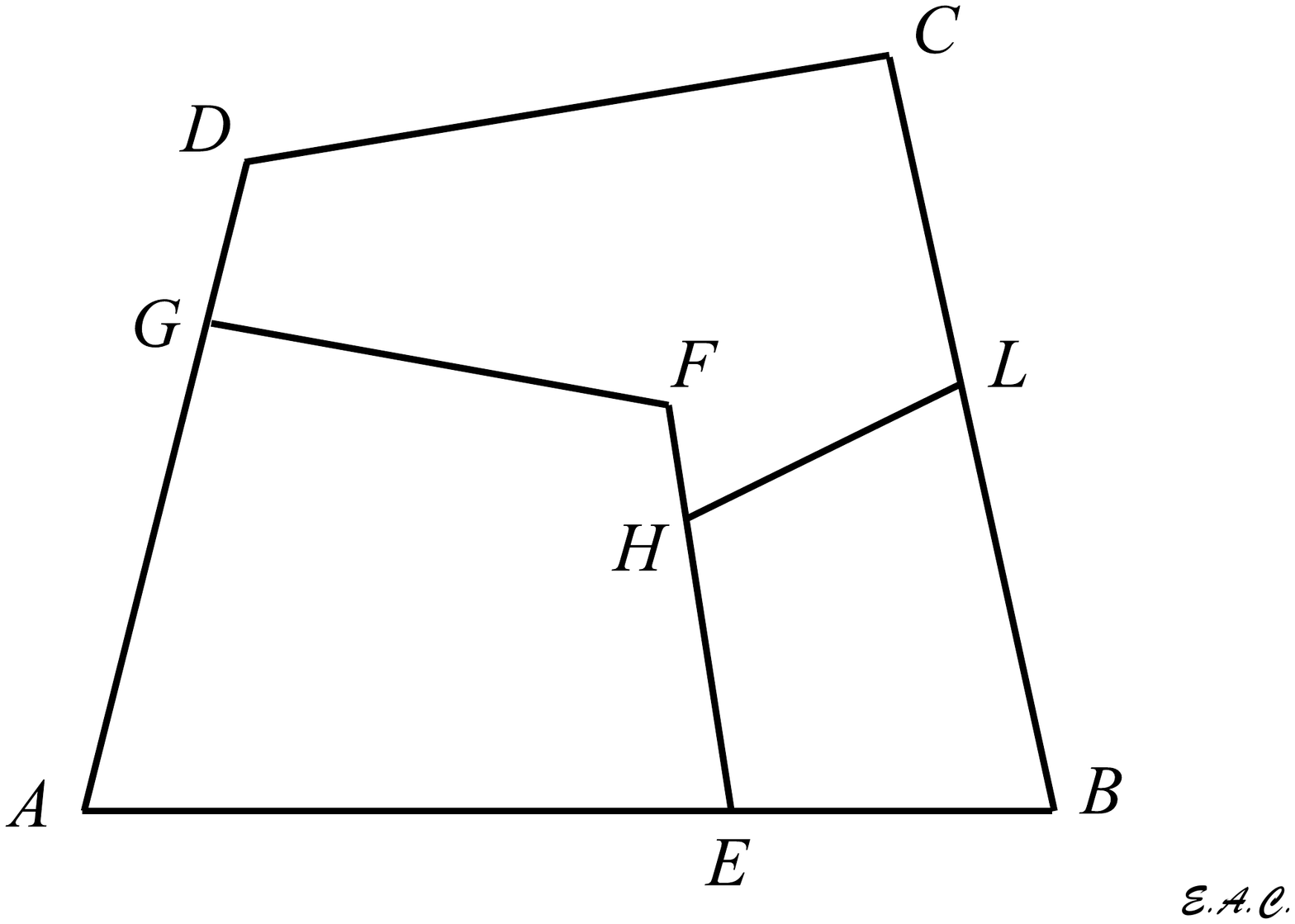,width=6.0in}}
\caption{}
\label{fig:planar}
\end{figure}
\section{(The planar articulated quadrilateral)}
It is not without interest to note the analogy between the theory of the
articulated octahedron with that of systems of articulated quadrilaterals 
studied by Mesrs. Hart (\cite{Hart}) and Kempe (\cite{Kempe})
for particular cases and Mr. Darboux (\cite{Darboux}) for the
general case. One of the problems studied by these geometers is, in effect, 
the following: given three planar quadrilaterals $ABCD$, $AEFG$, $BEHI$ 
arranged as shown in fig. \ref{fig:planar}, under what conditions does
their system, generally rigid, become susceptible to deformation?
If the tangents of the halves of the angles $\angle BAD$, $\angle ABC$, 
$\angle BEH$ are denoted by
$t, u, v$, there exist among these quantities three relations that have the 
same form as relations (\ref{tetra}), (\ref{tetrb}), (\ref{tetrc}) and they
ought to have an infinity of solutions. It is the same problem that presents
itself in the theory of the articulated octahedron \cite{note5}.

\appendix

\section{Comments on the preceding M\'emoire}
\begin{center}
{\bf by Mr. Am\'ed\'ee Mannheim.}
\end{center}
The M\'emoire of Mr. Bricard, very interesting in itself, is of particular
value from the point of view of {\em Kinematic Geometry}, because it
constitutes a chapter in the study of the displacement of a triangle
in space.

Mr. Bricard, by discovering the deformable octahedra, proved thus that 
under certain conditions a triangle of fixed size may be displaced in such 
a manner that its vertices describe circular arcs.

In effect, if we assume one of the facets $ABC$ of one of the deformable
octahedra to remain fixed, during the deformation the vertices $D$, $E$, $F$
of the opposite facet revolve as a result around the edges of the facet $ABC$
and the triangle $\triangle DEF$ is displaced. 
In the case of an arbitrary octahedron,
with $ABC$ fixed, the displacement of the triangle 
$\triangle DEF$ is not possible.
This is understood easily by noticing that the vertices of that triangle
must each describe circles, each of which is constrained by two conditions.
The triangle itself must then be constrained by six conditions: it is 
therefore immobile.

Let us consider (\ref{fig:triangle1}) a deformable octahedron whose facet 
$ABC$ is fixed. The triangle $\triangle DEF$ may therefore be displaced. 
The well
known properties, related to its displacement, lead to the properties of the
octahedron. I will give a single example.

Let us apply this theorem:
\vspace{0.1in}
\newline
{\it The normal planes to the trajectories of points of a plane pass through
a point of that plane.}
\vspace{0.1in}

The plane perpendicular to the trajectory of the point $D$ is the plane of the
facet $BDC$. Similarly, the planes of the facets that pass through $AB$, $AC$
are the normal planes to the trajectories of the vertices $F$, $E$. These
three planes meet therefore at a point of the plane $DEF$.

According to this, the following theorem can be stated:
\vspace{0.1in}
\newline
{\it The planes of the facets of a deformable octahedron, which pass through
the sides of one of the facets of that polyhedron, meet at a point of the
plane of the opposite facet.}
\vspace{0.1in}

This way, the properties of deformable octahedra may be found, but independently
of these properties, there appear additional questions related to a mobile
triangle.

For an infinitesimally small displacement of the triangle, there exists an
axis of displacement: what is the locus of these axes if the displacement
of the triangle is continued? What is the locus of the envelopes of the
plane of the mobile triangle? etc., etc.

It is seen that following the discovery of deformable octahedra, the study of
displacement of a triangle is presented with special conditions and
it requires new investigations by geometers.

\begin{thebibliography}{99}
\bibitem{Mannheim} These interesting comments were communicated to me by
Mr. Mannheim.
\bibitem{typo} (In the original, the first argument of $\phi$ in the 
equation above was omitted; it is clear that it expresses the relation
$v = \phi(t,u)$.)
\bibitem{note5} The models of the three octahedra, constructed from cardboard
leaves, have been offered to the Collections of the \'Ecole Polytechnique.
\end{thebibliography}
\end{document}